%% file: NoEigOutsideBulk6.tex
\documentclass[a4paper,12pt,reqno]{amsart}

\usepackage{etex}
\usepackage{pgfplots}
\usepackage{rotating,pdflscape,color,amssymb,subfigure,psfrag,amsmath,eufrak,bbm,epsfig}
\usepackage{a4wide,amscd}
\usepackage{tikz} 
 \usepackage[usenames,dvipsnames]{pstricks}
 \usepackage{pst-grad} 
 \usepackage{pst-plot} 

\definecolor{vdarkred}{rgb}{0.6,0,0.2}
\definecolor{vdarkblue}{rgb}{0,0.2,0.6}
\usepackage[pdftex, colorlinks, linkcolor=vdarkblue,citecolor=vdarkred]{hyperref}
\usepackage[small]{caption}

\newcommand{\me}{\mathrm{e}}
\newcommand{\ii}{\mathrm{i}}

\newcommand{\Cc}{\mc{C}}

\newcommand{\F}{\mc{F}}

\newcommand{\ld}{\ldots}
\newcommand{\beg}{\begin}
\newcommand{\en}{\end}

\newcommand{\trm}{\textrm}
 
\newcommand{\bgt}{\begin{itemize}}
\newcommand{\ent}{\end{itemize}}
\newcommand{\ite}{\item}

\newcommand{\op}{\operatorname}
\newcommand{\eqre}{\eqref}
\newcommand{\re}{\ref}
\newcommand{\la}{\label}

\newcommand{\diag}{\operatorname{diag}}

 \newcommand{\bgn}{\begin{enumerate}}
\newcommand{\enn}{\end{enumerate}}
\newcommand{\ds}{\displaystyle}
 
\newcommand{\p}{\mathbb{P}}

\newcommand{\tr}{ \operatorname{tr}}
\newcommand{\Tr}{\operatorname{Tr}}

\newcommand{\E}{\mathbb{E}}

\newcommand{\R}{\mathbb{R}}
\newcommand{\C}{\mathbb{C}}

\newcommand{\ud}{\mathrm{d}}

\newcommand{\pro}{probability }

\newcommand{\f}{\frac}
\newcommand{\ff}{\frac{1}}
\newcommand{\lf}{\left}
\newcommand{\ri}{\right}

\newcommand{\st}{such that }

\newcommand{\lam}{\lambda}

\newcommand{\ti}{\times}

\newcommand{\vfi}{\varphi}
\newcommand{\ste}{\, ;\, }
\newcommand{\mc}{\mathcal }

\newcommand{\eps}{\varepsilon}

\newcommand{\D}{\mc{D}}

\newcommand{\bck}{\backslash}
\newcommand{\al}{\alpha}

\newcommand{\eqlaw}{\stackrel{\textrm{law}}{=}}

\newcommand{\bbm}{\begin{bmatrix}}
\newcommand{\ebm}{\end{bmatrix}}
\newcommand{\bes}{\begin{equation*}}
\newcommand{\ees}{\end{equation*}}
\newcommand{\be}{\begin{equation}}
\newcommand{\ee}{\end{equation}}
\newcommand{\beqy}{\begin{eqnarray}}
\newcommand{\eeqy}{\end{eqnarray}}
\newcommand{\beq}{\begin{eqnarray*}}
\newcommand{\eeq}{\end{eqnarray*}}
\newcommand{\one}{\mathbbm{1}}
\newcommand{\lto}{\longrightarrow}

\newcommand{\ie}{\emph{i.e.\@ }}
\newcommand{\eg}{\emph{e.g.\@ }}
\newcommand{\bpm}{\begin{pmatrix}}
\newcommand{\epm}{\end{pmatrix}}

\newcommand{\cd}{\cdots}

\newcommand{\bpr}{\beg{proof}}
\newcommand{\epr}{\en{proof}}

\newcommand{\del}{\delta}

\newcommand{\pa}{\partial}

\newcommand{\tG}{\tilde{G}}

\newcommand{\ka}{\kappa}

\newcommand{\tg}{\tilde{g}}

\newcommand{\tQ}{\tilde{Q}}
\newcommand{\btQ}{\bar{\tilde{Q}}}

\newcommand{\oo}[1]{\overset{\circ}{#1}} 
\newcommand{\on}{\|\,\cdot\,\|} 
\newcommand{\Lc}{\mc{L}} 
\newcommand{\trans}{{\sf T}}

\newtheorem{Th}{Theorem}[]

\newtheorem{propo}[Th]{Proposition}

\newtheorem{lem}[Th]{Lemma}

\theoremstyle{definition}

\newtheorem{rem}[Th]{Remark}

 \allowdisplaybreaks
 
\date{\today}
\thanks{Couillet's work is supported by the ANR Project RMT4GRAPH (ANR-14-CE28-0006).}

\title[]{Spectral analysis of the Gram matrix\\
of mixture models}

\author{Florent Benaych-Georges}

\address[F.B.G.]{MAP 5, UMR CNRS 8145 -- Universit\'e Paris Descartes, Paris, France.}
\email{florent.benaych-georges@parisdescartes.fr}
\author{Romain Couillet}
\address[R.C.]{CentraleSup\'elec -- LSS -- Universit\'e ParisSud, Gif sur Yvette, France}
\email{romain.couillet@supelec.fr}



\begin{document}
\maketitle

\begin{abstract}This text is devoted to the asymptotic  study  of some spectral properties of the Gram matrix $W^{\sf T} W$ built upon a collection $w_1, \ldots, w_n\in \mathbb{R}^p$ of random vectors (the columns of $W$), as both the number $n$ of observations and the dimension $p$ of the observations tend to infinity and are of similar order of magnitude.  The random vectors $w_1, \ldots, w_n$ are independent  observations, each of them belonging to one of $k$ classes $\mathcal{C}_1,\ldots, \mathcal{C}_k$. The observations of  each class $\mathcal{C}_a$ ($1\le a\le k$)  are characterized by their distribution $\mathcal{N}(0,  p^{-1}C_a)$, where $C_1, \ldots, C_k$ are some non negative definite $p\times p$ matrices. The cardinality $n_a$ of class $\mathcal{C}_a$ and the dimension $p$ of the observations are    such that $n_a/n$ ($1\le a\le k$) and $p/n$ stay bounded away from $0$ and $+\infty$. We provide deterministic equivalents to the empirical spectral distribution of $W^\trans W$ and to the matrix entries of its resolvent (as well as of the resolvent of $WW^{\sf T}$). These deterministic equivalents are defined thanks to the solutions of a fixed-point system. Besides, we prove that $W^{\sf T} W$ has asymptotically no eigenvalues outside the bulk of its spectrum, defined thanks to these deterministic equivalents.  These results are directly used  in our companion paper \cite{FloRomainStat}, which is  devoted to the analysis of the spectral clustering algorithm in large dimensions. They also find applications in various other  fields    such as wireless communications where functionals of the aforementioned resolvents allow one to assess the communication performance across multi-user multi-antenna channels.\end{abstract}

\section{Introduction and main results}

\subsection{Introduction} In this article, we consider the Gram matrix $W^\trans W$, where $W=[w_1, \ld, w_n]$, $w_i\in \R^p$, is a collection of independent random vectors. Each $w_i$ belongs to one of $k$ classes $\Cc_1 ,\ld,\Cc_k$, with $\Cc_a$ ($1\le a\le k$) the class of vectors distributed as $\mc{N}(0,  p^{-1/2}C_a)$, where $C_1, \ld, C_k$ are some non negative definite $p\ti p$ matrices.  

This $k$-fold class setting comes in naturally in the field of machine learning and in statistical problems such as kernel spectral clustering (see the companion paper \cite{FloRomainStat}, where the present    results   find direct applications). Clustering algorithms are methods used to discover unknown subgroups or clusters in data: they seek partitions of the data into distinct groups so that the observations within each group are the most similar to one another and the observations belonging to different groups are the most dissimilar from each other. The most commonly used  such algorithms are $k$-means, the hierarchical clustering and EM \cite{HTF,JWTT}. Spectral clustering techniques (see \eg \cite{Luxburg,KannanVempala}) make use of the spectrum of a similarity matrix of the data or of a more involved derived version of it (such as the associated \emph{Laplacian matrix}) to perform dimensionality reduction before clustering in fewer dimensions, usually thanks to one of the previously mentioned algorithms. More specifically, given $n$ observations $x_1, \ld, x_n\in \R^p$ one wants to cluster, one chooses a similarity measurement, such as $\ka(x_i,x_j):=f(\frac{\|x_j-x_i\|^2}{p})$, for a well chosen decreasing function $f$, and defines the \emph{similarity matrix} $A:=\bbm\ka(x_i,x_j)\ebm_{i,j=1}^n$. Then, each observation  $x_i\in \R^p$ is replaced by its projection $y_i$ onto the linear span of the $k$ leading eigenvectors of the similarity matrix $A$ or of its   ``Laplacian" $L:=\bbm\f{\ka(x_i,x_j)}{\sqrt{d_id_j}}\ebm_{i,j=1}^n$,  where, for each $i$, $d_i:=\sum_j\ka(x_i,x_j)$, and $k$-means (or another   aforementioned standard algorithms) is then  performed on $y_1, \ld , y_n$. It appears in our companion article \cite{FloRomainStat} that, when the observations $x_i$ belong to classes, such that $x_i=\mu_a+p^{1/2}w_i$ for some $w_i\sim\mc{N}(0,p^{-1}C_a)$ when falling in class, say, $\Cc'_a$, the performance of the clustering algorithm relies on theoretical results on the spectrum and the resolvent of the aforementioned matrix $W^\trans W$, some of which are stated and proved here. 

Another class of applications is found in the field of wireless communications, where, letting $W_a=[w_{\sum_{j=1}^{a-1}n_j+1},\ldots,w_{\sum_{j=1}^an_j}]\in\R^{p\times n_j}$, $w_i\sim\mc{N}(0,p^{-1}C_a)$ be the communication channel between a $p$-antenna transmitter and an $n_a$-antenna receiver and $\sigma^2$ some ambient noise variance, functionals of the type $\log\det(WW^\trans+\sigma^2I_p)$ or $\tr W_a^\trans W_a(WW^\trans+\sigma^2I_p)^{-1}$ are instrumental to evaluate the maximally achievable communication rate across the channel, see \eg \cite{KAM09,COU11}.

The purpose of the present paper is to prove several technical results concerning the resolvent matrices of $W^\trans W$ and $WW^\trans$, which are then used, along with other arguments, to obtain a deterministic equivalent for the empirical spectral measure of $W^\trans W$ as $n,p\to\infty$ while $n/p$ remains away from $0$ and $\infty$, and to show that no eigenvalue can be found at macroscopic distance from the support of this measure. 
Some of these results and related questions (which we sometimes rederive for completion) can be found scattered in the literature: as in \cite{COU11}, where the same model is considered but only linear functionals of the empirical spectral measure of $W^\trans W$ are considered, or in \cite{MP67,CHO95,CouilletHachemRMTA,AEK15} where the spectral analyses of different models are considered, leading in particular to analogous equations to those introduced in Proposition~\re{defg_as} below (some of the results of \cite{MP67,CHO95} are immediate corollaries of our present results), or else in \cite{BaiSilver98,johnstone2001,BBP,BEN1,BEN4,CCHM,LV2011}, where the question of isolated eigenvalues of related matrix models was also considered.

 \subsection{Model presentation}
Let us start by introducing some notations. \emph{In what follows, everything, when not mentioned as \emph{constant} or \emph{fixed}, depends implicitly on the parameters $n,p,\ld$ introduced below.}

 Let $k$ be fixed, $n,p,n_1, \ld, n_k$ be some positive integers all tending to infinity \st $n_1+\cd+n_k=n$ and \st the ratios $$c_0=\f{p}{n},\qquad c_a=\f{n_a}{n}\quad (a=1, \ld, k)$$   all stay bounded away from $0$ and $+\infty$. We also introduce some positive semi-definite  $p\ti p$ matrices $C_1, \ld, C_k$ which we suppose to be bounded, in operator norm,  by a constant,  and some standard real Gaussian matrices $$Z_1\in \R^{p\ti n_1},\ld, 
 Z_k\in \R^{p\ti n_k}$$ and define the $p\ti n$ random  matrix $W$ by $$W:=p^{-1/2}[C_1^{1/2}Z_1,\ld, C_k^{1/2}Z_k]
.$$
In what follows, we will be interested in the \emph{empirical spectral distribution} $\mu_{W^\trans W}$ of $W^\trans W$ (\ie the uniform \pro measure on the spectrum of $W^\trans W$) and in the matrix entries of 
the \emph{resolvents}  $$Q_{z} :=(W^\trans W-z)^{-1}\qquad\trm{ and }\qquad \tQ_{z} :=(WW^\trans -z)^{-1}.$$

\subsection{Main results}

\input{Figure_Histogram.tex}

 Recall that everything, when not mentioned as \emph{constant}, depends implicitly on the parameters $n,p,z,\ld$ On the contrary, a \emph{constant} element is non random and independent of all other parameters.  Two elements are said to be at a \emph{macroscopic distance} from each other  if there is a constant $\eps>0$ \st there distance is lower-bounded by $\eps$.  We define the \emph{Stieltjes transform} of a \pro measure $\mu$ as the function $\ds m_\mu(z):=\int\f{ \mu(\ud t)}{t-z}.$ The set of signed measures on $\R$ is endowed with the \emph{weak topology}, \ie the topology defined by continuous bounded functions. 

\beg{Th}\la{Th:mu}  The measure $\mu$ defined by its Stieltjes transform   \be\la{defm_mu}m_\mu(z)=c_0\sum_{a=1}^kc_ag_a(z) ,\ee where the vector $(g_1(z), \ld, g_k(z))$, $z\in\C\setminus\R$, is defined by Proposition~\re{defg_as}, is a deterministic \pro measure   (depending on $n,p,\ld$) with compact support $\mathcal S$ \st   we have the almost sure convergences
\begin{align}
	\mu_{W^\trans W}-\mu &\lto 0 \label{eq:weak_conv} \\
\la{181513h29FOC}
	\operatorname{dist}\left(\operatorname{Spec}(W^\trans W), \mathcal S\cup \{0\} \right) &\lto 0.
\end{align}

Besides, if the matrices $C_1, \ld, C_k, I_p$ are all positive definite and linearly independent, $m_\mu(z)$ extends continuously from $\C^+:=\{z\in\C,\Im(z)>0\}$ to $\C^+\cup\R^*$ and, except for a possible atom at zero 
, $\mu$ has a continuous density given by $\ff{\pi}\Im(m_\mu(x))$ ($x\in\R^*$).
\en{Th}

\begin{rem}
	We believe that the technical assumption that $C_1, \ld, C_k, I_p$ are all positive definite and linearly independent to obtain a continuous limiting density can be relaxed. As it stands though, since an $\varepsilon$ perturbation of $C_1, \ld, C_k$ can always ensure these conditions while modifying $\operatorname{Spec}(W^\trans W)$ in a controlled manner, we have that for all large $n,p$, $\operatorname{Spec}(W^\trans W)$ is close to a spectrum having continuous deterministic equivalent and still defined by the equations of Proposition \re{defg_as} below, with some slightly perturbed versions of the matrices $C_a$. 
\end{rem}

\medskip

The vector $(g_1(z),\ldots,g_k(z))$ mentioned in Theorem~\ref{Th:mu} is defined explicitly in the following proposition.

\beg{propo}[Definition of $g_1(z),\ldots,g_k(z)$]\la{defg_as}For any $z\in \C\bck\R$, there is a unique vector  $(g_1(z), \ld, g_k(z))\in \C^k$ (depending on $n,p,\ld$) \st for each $a$,  \be\la{19152defpropo} \Im z\Im g_a(z)\ge  0,\qquad  \Im z\Im(zg_a(z))\ge  0,\qquad  c_0|g_a(z)|\le   (\Im z)^{-1}\ee and 
\begin{align}
	c_0g_a(z)  = -\frac1{z} \frac1{1+\tilde{g}_a(z)},\quad \tilde{g}_a(z) = -\frac1z \frac1p\operatorname{tr} C_a\left(I_p + \sum_{b=1}^k c_b g_b(z) C_b \right)^{-1}.
	\label{defg_adefpropo}
\end{align}
Besides, the first two inequalities in \eqre{19152defpropo} are in fact strict and  the functions  $c_0g_a(z)$ ($a=1, \ld, k$) are the Stieltjes transforms of some $\R_+$-compactly supported \pro  measures $\nu_1, \ld, \nu_k$.
\en{propo} 

\beg{rem}\la{numerics_for_mu}It appears in the proof of Proposition \re{defg_as} that the functions $g_a$ (hence the function $m_\mu$, by \eqre{defm_mu}) can be computed numerically in a  very efficient way thanks to a fixed point convergence.   Then, choosing $z=x+\mathrm{i}\eta$ close to the real line, we get a good approximation of the measure $\mu$ of Theorem \re{Th:mu} as $\mu\approx \ff{\pi}\Im(m_\mu(x+\mathrm{i}\eta))\ud x$.
\en{rem}

\medskip

At the core of the proof of Theorem~\ref{Th:mu} is the following result.
For $A=A(n),B=B(n)$ some random square matrices with     size $n$  tending to infinity, the notation $A \leftrightarrow B  $ stands for the fact that  we have the convergences in \pro $\frac1n\tr D (A - B) \lto 0$ and $d_{1 }^\trans (A -B )d_{2 }\lto 0$ for all sequence $D=D(n)$ of deterministic   $n\ti n$ matrices   of bounded norms and all deterministic sequnces  vectors $d_i=d_{i }(n)$ of bounded norms. 

\beg{propo}[Deterministic Equivalents]\la{prop:288151}Let  $\mc{S}:=\cup_{a=1}^k\op{supp}(\nu_a)$. For any $z\in\C$   at macroscopic distance from $ \mc{S}\cup\{0\}$, we have  
	\begin{align}\la{181513h29}
		Q_{z}  & \leftrightarrow \bar{Q}_{z}  := c_0 \diag \left\{ g_a(z) 1_{n_a} \right\}_{a=1}^k  \\
	\la{181513h292}	\tilde{Q}_{z} & \leftrightarrow \bar{\tilde{Q}}_{z}  := -\frac1z \left(    I_p + \sum_{a=1}^k c_a g_a(z) C_a  \right)^{-1}.
	\end{align}
\en{propo}

Proposition~\ref{prop:288151} by itself finds immediate applications in the aforementioned area of wireless communications, where functionals of the type $-\frac1p\log Q_x$ or $w_i^*Q_xw_i$ are fundamental quantities to evaluate achievable communication rates \cite{COU11}, or in machine learning where some more involved linear statistics of the Gram matrix $XX^\trans$, with $X=p^{-\frac12}[\mu_11_{n_1}^\trans,\ldots,\mu_k1_{n_k}^\trans]+W$ a Gaussian mixture sample, carries relevant information for data classification and clustering \cite{FloRomainStat}. For these applications, further results such as central limit theorems \cite{KAM09} or higher order statistics are required. For completion, we provide below some additional results that come in handy into this scope.

 \begin{propo}[Further Deterministic Equivalents]
	\label{prop:deteq2}
	For any $z_1,z_2\in\C$ at macroscopic distance from $ \mc{S}\cup\{0\}$, 
	\begin{align*}
		Q_{z_1}\mathcal D_a Q_{z_2} &\leftrightarrow \bar{Q}_{z_1}\mathcal D_a \bar{Q}_{z_2} + \sum_{b=1}^k R_{ab}(z_1,z_2) \bar{Q}_{z_1}  \mathcal D_b   \bar{Q}_{z_2} \\
		\tilde{Q}_{z_1}C_a \tilde{Q}_{z_2} &\leftrightarrow \bar{\tilde Q}_{z_1}C_a \bar{\tilde Q}_{z_2} + \sum_{b=1}^k R_{ba}(z_1,z_2) \bar{\tilde Q}_{z_1} C_b \bar{\tilde Q}_{z_2} \\
		\frac1p \tilde{Q}_{z_1} W\mathcal D_a W^\trans \tilde{Q}_{z_2} &\leftrightarrow z_1z_2 c_0 c_a g_a(z_1)g_a(z_2) \bar{\tilde Q}_{z_1}C_a \bar{\tilde Q}_{z_2}
	\end{align*}
	where $R(z_1,z_2)_{ab}=\frac{c_a}{c_b} [(I_k-\Omega(z_1,z_2))^{-1}\Omega(z_1,z_2)]_{ab}$ with, for $1\leq a,b\leq k$,
	\begin{align*}
		\Omega(z_1,z_2)_{ab} &=  c_0 c_b z_1 g_a(z_1) z_2g_a(z_2) \frac1p\tr C_a \bar{\tilde Q}_{z_1} C_b \bar{\tilde Q}_{z_2}.
	\end{align*}
\end{propo}

Although it appears only at this point, the matrix $\Omega(z_1,z_2)$ is a cornerstone of the proof of Theorem~\ref{Th:mu}. It is in particular related to the derivative of the $g_a(z)$'s introduced in Proposition~\ref{defg_as} as follows
\begin{align*}
	\left\{ g_a'(z) \right\}_{a=1}^k &= c_0 \left( I_k - \Omega(z,z)\right)^{-1} \left\{ g_a^2(z) \right\}_{a=1}^k.
\end{align*}

 \section{Proofs: preliminary approximation lemma}
 Here, we denote the Hadamard product (\ie entry-wise product)  of matrices or vectors by $\odot$. We shall also use both superscript notations $(\cdot)^\trans$ and $(\cdot)^*$ to denote transpose and Hermitian transpose of matrices and vectors.

 \beg{lem}\la{lem:keyapproxlemma}Let $z,\underline{z}\in \C\bck\R$ and   $(g_a)_{a=1}^k, (\underline{g}_a)_{a=1}^k\in \C^k$ \st for each $a$, \be\la{19152}\Im(z)\Im(g_a)\ge 0,\qquad \Im(z)\Im(zg_a)\ge 0,\qquad  c_0|g_a|\le |\Im z|^{-1}\ee(and the same for $\underline{g}_a,\underline{z}$)   and $(\eps_a)_{a=1}^k$, $(\underline{\eps}_a)_{a=1}^k$, $(\eta_a)_{a=1}^k$, $(\underline{\eta}_a)_{a=1}^k\in \C^k$ \st for each $a=1,\ld, k$, we have $$c_0g_a  = \ff{-z (1+\ff{p}\tr C_a\btQ+\eta_a)}+c_0\eps_a \quad; \quad  c_0\underline{g}_a  = \ff{-\underline{z} (1+\ff{p}\tr C_a\underline{\btQ}+\underline{\eta}_a)}+c_0\underline{\eps}_a,$$ with $$\btQ:=-z^{-1}\lf(I_p+\sum_{a=1}^kc_ag_aC_a\ri)^{-1}\qquad ;\qquad  \underline{\btQ}:=-{\underline{z}}^{-1}\lf(I_p+\sum_{a=1}^kc_a\underline{g}_aC_a\ri)^{-1}.$$
 Then:
 \bgt\ite[(i)] We have $$(I-\Omega)(\underline{g}-g)=c_0(\underline{z}-z+\underline{z}\underline{\eta}-z\eta)\odot(\underline{g}-\underline{\eps})\odot (g-\eps)+\underline{\eps}-\eps$$ for $$\Omega=\Omega(z,\underline{z},g,\underline{g},\eps,\underline{\eps},\eta,\underline{\eta}):=[c_0z\underline{z}(\underline{g}_a-\underline{\eps}_a)(g_a-\eps_a)c_b\ff{p}\tr C_a\btQ C_b\underline{\btQ}]_{a,b=1}^k.$$ 
 \ite[(ii)] For $C:=\max _a \|C_a\| $, 
 we have $$\rho(\Omega)\le 1-\min\lf\{ \f{(\Im z)^2}{|z|( |\Im z|+C)} ,\f{(\Im \underline{z})^2}{|\underline{z}|( |\Im \underline{z}|+C)} \ri\}^2+\|\eps\|_\infty+\|\underline{\eps}\|_\infty+c_0(|z|\|\eta\|_\infty+|\underline{z}|\|\underline{\eta}\|_\infty).$$
 \ite[(iii)]  There are $P,Q$ some polynomials with non negative coefficients and  $c_1>0$   \st for  $\al:=\max\{\|\eps\|_\infty, \|\underline{\eps}\|_\infty,\|\eta\|_\infty,\|\underline{\eta}\|_\infty\}$, we have 
 $$\al\le c_1 \min\{|\Im z|^{4},|\Im \underline{z}|^{4}, |z|^{-3},|\underline{z}|^{-3}\} \implies \|\underline{g}-g\|\le (|\underline{z}-z|+\al)P(|z|+|\underline{z}|)Q(|\Im z|^{-1}+|\Im \underline{z}|^{-1}).$$
 \ent
 \en{lem}

 \bpr  
 Note first that by the hypotheses $zg_a,\underline{z}\underline{g}_a\in \C^+$, we have that $\|\btQ\|\le |\Im z|^{-1}$. It follows that \be\la{19151}|c_0g_a|\ge |z|^{-1}\ff{1+C|\Im z|^{-1}+|\eta_a|},\ee for $C=\max _a \|C_a\|_\infty$. The same kind of inequalities hold for $\underline{\btQ}$ and $\underline{g}_a$.

We have 
\beq c_0(\underline{g}_a-g_a)
&=&c_0(\underline{\eps}_a-\eps_a)+ c_0(\underline{g}_a-\underline{\eps}_a)-c_0(g_a-\eps_a)\\
&=&c_0(\underline{\eps}_a-\eps_a)+  \ff{z (1+\ff{p}\tr C_a\btQ+\eta_a)}- \ff{\underline{z} (1+\ff{p}\tr C_a\underline{\btQ}+\underline{\eta}_a)}\\
&=&c_0(\underline{\eps}_a-\eps_a)+ \f{\underline{z}-z+\underline{z}\underline{\eta}_a-z\eta_a+\ff{p}\tr (C_a(\underline{z}\underline{\btQ}-z\btQ))}{z\underline{z} (1+\ff{p}\tr C_a\btQ+\eta_a)  (1+\ff{p}\tr C_a\underline{\btQ}+\underline{\eta}_a)}.
\eeq
Now, one has to notice that as  $\underline{\btQ}^{-1}+\underline{z}I_p=-\underline{z}\sum_{b=1}^kc_b\underline{g}_bC_b$, we have 
  \beq z\tr C_a\btQ+z\underline{z}\tr C_a\btQ\underline{\btQ}&=& z\tr C_a\btQ(\underline{\btQ}^{-1}+\underline{z})\underline{\btQ}\\
  &=&-z\underline{z}\sum_{b=1}^kc_b\underline{g}_b \tr C_a\btQ C_b\underline{\btQ}
  \eeq   
  In the same way, as  $\btQ^{-1}+zI_p=-z\sum_{b=1}^kc_bg_bC_b$, we have 
  \beq \underline{z}\tr C_a\underline{\btQ} +z\underline{z}\tr C_a\btQ\underline{\btQ}&=& z\tr C_a\btQ(\btQ^{-1}+z)\underline{\btQ}\\
  &=&-z\underline{z}\sum_{b=1}^kc_bg_b\tr C_a\btQ C_b\underline{\btQ}
  \eeq  
It follows that 
\beq \tr (C_a(\underline{z}\underline{\btQ}-z\btQ))&=& \underline{z}\tr C_a\underline{\btQ} +z\underline{z}\tr C_a\btQ\underline{\btQ}-(z\tr C_a\btQ+z\underline{z}\tr C_a\btQ\underline{\btQ})\\
&=&z\underline{z}\sum_{b=1}^kc_b\underline{g}_b \tr C_a\btQ C_b\underline{\btQ}-z\underline{z}\sum_{b=1}^kc_bg_b\tr C_a\btQ C_b\underline{\btQ}
\\
&=&z\underline{z}\sum_{b=1}^kc_b(\underline{g}_b-g_b)\tr C_a\btQ C_b\underline{\btQ}
\eeq
and that 
\beq c_0(\underline{g}_a-g_a)
&=&c_0(\underline{\eps}_a-\eps_a)+ \f{\underline{z}-z+\underline{z}\underline{\eta}_a-z\eta_a+\ff{p}\tr (C_a(\underline{z}\underline{\btQ}-z\btQ))}{z\underline{z} (1+\ff{p}\tr C_a\btQ+\eta_a)  (1+\ff{p}\tr C_a\underline{\btQ}+\underline{\eta}_a)}\\
&=&c_0(\underline{\eps}_a-\eps_a)+c_0^2(\underline{g}_a-\underline{\eps}_a)(g_a-\eps_a)\times \\ &&\lf(\underline{z}\underline{\eta}_a-z\eta_a+(\underline{z}-z)+z\underline{z}\sum_{b=1}^kc_b(\underline{g}_b-g_b)\ff{p}\tr C_a\btQ C_b\underline{\btQ}\ri)
\eeq
We directly deduce (i).

To prove (ii), let us first treat the case  where  $\underline{z}=z^*$ and $\underline{\eps}_a=\eps_a^*$, $ \underline{g}_a=g_a^*$ for each $a$. If $\Im(z)>0$ (the other case can be treated in the same way), it is easy to see that $\Omega$ has positive entries (as $\tr C_a\btQ C_b\underline{\btQ}=\tr C_a^{1/2}\btQ C_b^{1/2}(C_a^{1/2}\btQ C_b^{1/2})^*$) and that we have, by   (i),  $$(I-\Omega)\Im(g)=c_0\Im(z)|g-\eps|^2+\underbrace{\f{c_0}{2} (\underline{\eta}\underline{z}-\eta z)+\ff{2} (\underline{\eps}-\eps)}_{\ds :=\ka}.$$
Thus $\ka\in \R^k$ and  if $u=(u_1, \ld, u_k)$ is a left eigenvector of $\Omega$ associated with $\rho(\Omega)$ \st for all $a$, $u_a\ge 0$ and $\sum_{a=1}^k u_a=1$ (as per Lemma~\ref{lem:nonnegative_matrices} in the appendix),  then we have $$(1-\rho(\Omega))\sum_{a=1}^ku_a\Im (g_a)=c_0\Im(z)\sum_{a=1}^ku_a|g_a-\eps_a|^2+u^\trans \ka,$$ which implies, using successively \eqre{19152} and \eqre{19151}, \beq 1-\rho(\Omega)&=&c_0\f{\sum_{a=1}^ku_a|g_a|^2}{\sum_{a=1}^ku_a\Im (c_0g_a)(\Im(z))^{-1}}+u^\trans\ka\\ &\ge& c_0^2\f{\sum_{a=1}^ku_a|g_a|^2}{\sum_{a=1}^ku_a (\Im(z))^{-2}}+u^\trans\ka\\ &=&(c_0 \Im(z))^2\sum_{a=1}^ku_a|g_a|^2+u^\trans\ka\\ &\ge& \lf(\f{(\Im z)^2}{|z|( |\Im z|+C)}\ri)^2+u^\trans\ka,\eeq
 so that the spectrum of $\Omega$ is contained in the ball with center $0$ and radius $$1-\lf(\f{(\Im z)^2}{|z|( |\Im z|+C)}\ri)^2+\|\eps\|_\infty+c_0|z|\|\eta\|_\infty.$$ 
 
 To treat the general case, just use Lemmas \re{lem:spectral_radius} and \ref{lem:Cauchy_Schwarz_spectral_radius} from the appendix and notice that $$|\tr C_a\btQ C_b\underline{\btQ}|=|\tr C_a^{1/2}\btQ C_b^{1/2}C_b^{1/2}\underline{\btQ}C_a^{1/2}|\le \sqrt{\tr C_a\btQ C_b\btQ \tr C_a\underline{\btQ} C_b\underline{\btQ}}$$ from which it follows that $\rho(\Omega(z,\underline{z},\ldots))\leq \sqrt{\rho(\Omega(z,z^*,\ldots)) \rho(\Omega(\underline{z},\underline{z}^*,\ldots))}$ and then use $\sqrt{AB}\leq\max(A,B)$ for $A,B>0$.

At last, (iii) follows from the formula of the inverse of a matrix in terms of the determinant and of the minors. 
 \epr

\section{Proof of Proposition \re{defg_as}}

\subsection{Uniqueness}Note that for each fixed $z\in \C\bck\R$, if there exist two vectors $$(g_1(z), \ld, g_k(z))\qquad \trm{ and }\qquad (\underline{g}_1(z), \ld, \underline{g}_k(z))$$ satisfying \eqre{19152defpropo} and the equations \eqre{defg_adefpropo}, then one can apply Lemma \re{lem:keyapproxlemma} with $\eps=\eta=0$: by (i), we get that $(I-\Omega)(\underline{g}-g)=0$, whereas by (ii) we know that $\rho(\Omega)<1$, which implies that $g=\underline{g}$. 

\subsection{Existence}Note first that one can focus on $\C^+$ and then extend to $\C\bck\R$ by the formula $g_a(z^*)=g_a(z)^*$. We shall first prove that there is a unique collection of functions of $z$ satisfying the  weakened version of  conditions \eqre{19152defpropo} given by \eqre{19152def} and the equations \eqre{defg_adefpropo}.

Let $\mc{L}$ be the set of analytic functions $g:\C^+\to \C$ \st for all $z\in \C^+$, 
\be\la{19152def} \Im g(z)\ge 0,\qquad  \Im(zg(z))\ge 0,\qquad  c_0|g(z)|\le |\Im z|^{-1}\ee For   $\eta>0$,  we define, for $g,\underline{g}\in \Lc$, $$\op{d}_{\Lc,\eta} (g,\underline{g}):=  \sup_{\Im z\ge \eta}|\underline{g}(z)-g(z)| .$$ By   the analytic continuation principle, this is obviously a distance. 
 \beg{lem}$\op{d}_{\Lc,\eta}$ is a complete distance on $\Lc$.  
\en{lem}

\bpr     Note first that by Montel theorem $\Lc$ is a compact subset of the set of analytic functions on $\C^+$ endowed with the topology of uniform convergence on compact sets. Let $(g_n)$ be a Cauchy sequence in $(\Lc,\op{d}_{\Lc,\eta})$. Then there is an analytic function $g$ defined on $\{z\ste\Im z>\eta\}$ \st on $\{z\ste \Im z>\eta\}$, $(g_n)$ converges uniformly to $g$. Besides, any accumulation point of $(g_n)$ in the set of analytic functions on $\C^+$ endowed with the topology of uniform convergence on compact sets coincides with $g$ on $\{z\ste \Im z>\eta\}$, hence $g$ is the restriction to $\{z\ste\Im z>\eta\}$ of  an element  of $\Lc$, and $(g_n)$ converges to $g$ in $(\Lc,\op{d}_{\Lc,\eta})$.
\epr

 We denote $\Lc^k=\Lc\ti\cdots\ti \Lc$.

\beg{lem}\la{lemfixedPoint}Let $\Psi:\Lc^k\to\Lc^k$ be defined by  $\Psi(g_1, \ld, g_k)=(f_1, \ld, f_k),$ where \be\la{defg_adef}c_0f_a(z)  = -\frac{1}{z  -\frac{1}{p}\operatorname{tr} C_a(I_p + \sum_{b=1}^k c_b g_b(z) C_b)^{-1}} \qquad (a=1,\ld, k)
	\ee
 Then $\Psi$ is well defined and  admits a unique fixed point in $\Lc^k$. 
\en{lem}

\bpr Let us first make two remarks. First, for   $C,D$ some non negative definite Hermitian matrices, $\tr CD=\tr C^{1/2}D^{1/2}( C^{1/2}D^{1/2})^*\ge 0$. By linear combination, if $C,D$ are only Hermitian matrices, $\tr CD\in \R$. Secondly, let $A$ be an invertible matrix \st $A=X+\ii Y$, with $X,Y$ Hermitian matrices \st $Y$ non negative definite. Then  $$A^{-1}=A^{-1}A^*(A^{-1})^*=A^{-1}X(A^{-1})^*-\ii A^{-1}Y(A^{-1})^*$$ has a skew-Hermitian part which has the form $\ii$ times a  non positive definite matrix. 

From both of these remarks, we deduce that for $g_1, \ld, g_k\in \C,$ we have $$\Im g_1\ge 0, \ld, \Im g_k\ge 0\implies \Im \operatorname{tr}  C_a(I_p + \sum_{b=1}^k c_b g_b C_b)^{-1}\le 0$$ and $$\Im (zg_1)\ge 0, \ld, \Im(z g_k)\ge 0\implies \Im \operatorname{tr} C_a(zI_p + \sum_{b=1}^k c_b zg_b C_b)^{-1}\le 0,$$
so that $\Psi:\Lc^k\to\Lc^k$ is well defined.

Let now $\eps>0$ \st for any $p\ti p$ matrices $X,Y$, $$\|Y-I_p\|, \|X-I_p\|\le \eps\implies \|Y^{-1}-X^{-1}\|\le 2\|Y-X\|.$$ Let $\eta_0>0$ be such that $$\sum_{b=1}^kc_b\eta_0^{-1}\|C_b\|\le \eps.$$

Now, fix $g,\underline{g}\in \Lc^k$ and set $f,\underline{f}:=\Psi(g), \Psi(\underline{g})$. For any $z\in \C^+$ \st $\Im z\ge \eta_0$,  and any $a=1, \ld, k$, we have, for $G:= \sum_{b=1}^k c_b g_b(z) C_b$ and $\underline{G}:= \sum_{b=1}^k c_b g_b(z) C_b$,
\beq|\underline{f}_a(z)-f_a(z) |&=& \f{c_0^{-1}}{|(z-\ff{p}\tr \f{C_a}{I_p+G})(z-\ff{p}\tr \f{C_a}{I_p+\underline{G}})|} \lf|\ff{p}\tr \lf(C_a((I_p+G)^{-1}-(I_p+\underline{G})^{-1}) \ri)\ri|  \\
&\le & 2c_0^{-1}k(\Im z)^{-2}\max_a\|C_a\|^2\max_a|\underline g_a(z)-g_a(z)|
\eeq
We deduce that for $\eta_0$ large enough, if one endows $\Lc^k$ with the product distance defined by $\op{d}_{\Lc,\eta}$,   then $\Psi$ is a contraction. By the previous lemma, it allows one to conclude.
\epr

\subsection{Conclusion of the proof}The functions $g_1, \ld g_k$ are analytic on $\C^+$  and satisfy \eqre{19152def} and the equations \eqre{defg_adefpropo}. Using $\Im(zg(z))\ge 0$ in \eqre{defg_adefpropo}, one easily gets that for each $a$, $$\lim_{\Im z\to+\infty} c_0zg_a(z)=-1.$$ This proves that the $c_0g_a$'s are the Stieltjes transforms of some \pro measures $\nu_a$, hence that for each $a$, $$
\Im g_a(z)> 0,\qquad     c_0|g_a(z)|\le |\Im z|^{-1}$$ Besides, as $$ \Im(zg_a(z))\ge 0, $$ the $g_a$'s are supported on $\R_+$, hence $ \Im(zg_a(z))> 0.$

Now, it remains to prove that the $\nu_a$'s have compact supports. For any $d,\eps>0$, let $\F_{d,\eps}$ denote the set of continuous functions $g:\{z\ste\Re z\ge d, \Im z\ge 0\}\to \C$ \st for all $z\in \R$, $g(z)\in \R$ and  for all $z\notin \R$, $$\Im g(z)\ge 0,\qquad  \Im(zg(z))\ge 0,\qquad   |g(z)|\le  \eps.$$ Clearly, when endowed with the distance $$\op{d}_{\F_{d,\eps}}(g,\underline g):=\sup_{\Re z\ge d, \Im z\ge 0} |\underline g(z)-g(z)|,$$ $\F_{d,\eps}$ is a complete metric space. 
Let $\Phi: \F_{d,\eps}^k\to \F_{d,\eps}^k$  be defined by  $\Phi(f_1, \ld, f_k)=(h_1, \ld, h_k),$ where \be\la{defPhi}c_0h_a(z)  = -\frac{1}{z  -\frac{1}{p}\operatorname{tr} \frac{C_a}{I_p + \sum_{b=1}^k c_b f_b(z) C_b}} \qquad (a=1,\ld, k)
	\ee
 Then by the same kind of computations as in the proof of Lemma \re{lemfixedPoint}, one proves that for $a$ large enough and $\eps$ small enough,  $\Phi$ is well defined and  admits a unique fixed point in $\F_{d,\eps}^k$.
By the pointwise uniqueness we already proved, this fixed point must coincide with the vector $(g_1,\ld, g_k)$ on $\{z\ste\Re z\ge d, \Im z> 0\}$. We deduce that the functions $g_1,\ld, g_k$ can be extended as continuous functions on $\C^+\cup [d, +\infty)$ taking real values on $[d,+\infty)$. By \cite[Th. 2.4.3]{agz}, we conclude that the measures $\nu_1, \ld, \nu_k$ have supports contained in $[0,d]$.

\section{Proof of Proposition~\re{prop:288151}}
In this section we shall use the notation $u_p=O(v_p)$ for a sequence $u_p$ possibly depending on other parameters ($i,i',z,\ld$) \st there are some polynomials $P,Q$ with non negative coefficients \st uniformly in all parameters, $$u_p\le v_p P(|z|)Q(|\Im z|^{-1}).$$  
For $M=M_p$ a matrix, $M_p=O_{\on}(v_p)$ means that the operator norm $\|M\|$ of $M$ satisfies  $$\|M\|=O(v_p).$$
Also, for $X$ a (possibly multidimensional) random variable, we set $$\oo X:=X-\E X.$$
At last, for each $a=1, \ld, k$, we define $j_a\in \R^{n\ti 1}$ as the column vector with $i$-th entry equal to $1$ if $n_1+\cd+n_{a-1}<i\le n_1+\cd+n_a$ and to $0$ otherwise and set $\D_a$, the diagonal $n\ti n$ matrix with diagonal $j_a$.
 
\subsection{Boundedness of $WW^\trans$}
We shall use the following lemma, following from the fact that $$WW^\trans =p^{-1}\sum_{a=1}^kC_a^{1/2}Z_aZ_a^\trans C_a^{1/2}$$ and that as is well known (\emph{e.g.} by \cite[Lem. 7.3]{ESY2}), there are $t_0,c_0>0$ constant \st for all $t>0$, \be\la{18151}\p(\max_{a=1,\ld,k}p^{-1}\|Z_aZ_a^\trans\|>t_0+t)\le e^{c_0n(t-t_0)}.\ee

\beg{lem}\la{lem:boundWW}There are $t_0,c_0>0$ constant \st for all $t>0$, \be\la{181514h}\p(\max_{a=1,\ld,k}\|WW^\trans \|>t_0+t)\le e^{c_0n(t-t_0)}.\ee
\en{lem}

\subsection{Loop equations}\la{sec:loop_equations}
We shall prove the following lemma in the next sections. For short, we shall denote $Q_z$ by $Q$. 
\beg{lem}\la{lem:2881516h}The matrix $\E Q$ is a diagonal matrix with diagonal entries that are constant along the classes, \ie of the form $\sum_{a=1}^k \al_a\D_a$, with $\al_a\in \C^+$. Besides, for $z\in \C\bck\R$,
\beqy\la{181515h}\E Q&=&-z^{-1}(I_n+\E D)^{-1}+O_{\on}(p^{-1})\\
\la{181515h1} \E \tQ  
        &=& -z^{-1}\left(I_p+\sum_{a=1}^k\E \ff{p}\tr (\mc{D}(j_a)Q)C_a \right) ^{-1}  +O_{\on}(p^{-1})
\eeqy
where \be\la{238152}D:=\sum_{a=1}^k \ff{p}\tr (\tQ C_a)\D_a.\ee
\en{lem}

The fact that the matrix $\E Q$ is diagonal follows from the Neumann expansion for $|z|$ large enough (and from analytic continuation for small $z$) and from the fact that the $Z_a$'s are independent and with symmetric distribution. The fact that it is of the form given here follows from the invariance of the law of $W^\trans W$ under conjugation by the appropriate permutation matrices (those with all cycles contained in a class).

For each $i=1, \ld, n$, we denote by $C(i)$ the covariance matrix of the $i$-th column of $p^{1/2}W$, so that $C(i)=C_a$ if $n_1+\cd+n_{a-1}<i\le n_1+\cd+n_{a-1}+n_a$. 
\subsubsection{Computations on $Q$}
 
  By the resolvent identity, $$Q=z^{-1}W^\trans WQ-z^{-1}I_n$$ so that, using Stein Lemma (see in appendix),
 \beq \E Q_{ii'}&=&z^{-1}  (W^\trans WQ)_{ii'}-\del_{ii'}z^{-1}\\ &=&-\del_{ii'}z^{-1}+z^{-1}\E\sum_{j,l} w_{ji}w_{jl}Q_{li'}\\
 &=&-\del_{ii'}z^{-1}+z^{-1}\sum_{j,l,m}\E w_{ji}w_{mi} \E \pa_{w_{mi}}w_{jl}Q_{li'}\\
  &=&-\del_{ii'}z^{-1}+(pz)^{-1}\sum_{j,l,m}C(i)_{jm} \{ \one_{(m,i)=(j,l)}\E  Q_{li'}-\E w_{jl}(Q([\del_{ri}w_{ms}+\del_{si}w_{mr}]_{r,s=1}^n)Q)_{li'}\}
 \eeq
 where we used the fact that $\E w_{ji}w_{mi} =p^{-1}C(i)_{jm}$. With $E_{ij}$ the matrix with unique non-zero entry $[E_{ij}]_{ij}=1$, $$\pa_{w_{mi}}Q=-Q(\pa_{w_{mi}} W^\trans W)Q=-Q(E_{im}W+W^\trans E_{mi})Q=-Q([\del_{ri}w_{ms}+\del_{si}w_{mr}]_{r,s=1}^n)Q$$ 
from which, after replacement, we get \beq \E Q_{ii'}
  &=&-\del_{ii'}z^{-1}+(pz)^{-1}\sum_{j,l,m}C(i)_{jm} \{ \one_{(m,i)=(j,l)}\E  Q_{li'}-\E w_{jl}(Q([\del_{ri}w_{ms}+\del_{si}w_{mr}]_{r,s=1}^n)Q)_{li'}\}\\
        &=&-\del_{ii'}z^{-1}+(pz)^{-1}\E Q_{ii'}\tr C(i)-(pz)^{-1}\sum_{j,l,m,s}C(i)_{jm} \E w_{jl}Q_{li} w_{ms}Q_{si'} \\
        &&
        -(pz)^{-1}\sum_{j,l,m,r}C(i)_{jm} \E w_{jl}Q_{lr} w_{mr}Q_{ii'} \\
            &=&-\del_{ii'}z^{-1}+(pz)^{-1}\E Q_{ii'}\tr C(i)-(pz)^{-1}\sum_{j,m}C(i)_{jm} \E (WQ)_{ji} (WQ)_{mi'} \\
        &&
        -(pz)^{-1}\sum_{j,m}C(i)_{jm} \E (WQW^\trans )_{jm} Q_{ii'} \\
             &=&-\del_{ii'}z^{-1}+(pz)^{-1}\E Q_{ii'}\tr C(i)-(pz)^{-1}  \E (QW^\trans  C(i)WQ)_{ii'} \\
        &&
        -(pz)^{-1} \E \tr ( WQW^\trans  C(i)) Q_{ii'} 
        \eeq
   Besides, it is easy to see that $WQW^\trans =WW^\trans \tQ=z\tQ+I_p$ which entails
    \beq \E Q_{ii'}
  &=&-\del_{ii'}z^{-1}+(pz)^{-1}\E Q_{ii'}\tr C(i)-(pz)^{-1}  \E (QW^\trans  C(i)WQ)_{ii'} \\ &&
        -(pz)^{-1} \E \tr ( (z\tQ+1) C(i)) Q_{ii'}\\
          &=&-\del_{ii'}z^{-1}+(pz)^{-1}\E Q_{ii'}\tr C(i)-(pz)^{-1}  \E (QW^\trans  C(i)WQ)_{ii'} 
          -p^{-1}\E \tr (\tQ C(i)) Q_{ii'}
          \\ &&
        -(pz)^{-1} \tr (   C(i)) \E Q_{ii'}\\
            &=&-\del_{ii'}z^{-1}-(pz)^{-1}  \E (QW^\trans  C(i)WQ)_{ii'} 
	    -\E p^{-1}\tr (\tQ C(i)) Q_{ii'}
         \eeq
         In other words, if one defines \be\la{def:M1}M_1:= \ff{z}\sum_{a=1}^k \E \D_aQW^\trans  C_aWQ = \ff{z}\sum_{a=1}^k \E \D_aW^\trans \tQ C_a\tQ W,
         \ee then we have, for $D$ as in \eqre{238152} and $\oo{D}:=D-\E D$,  
    \beqy\la{218159h36} -z\E Q        &=& I_n+z\E DQ+zp^{-1}M_1
          \\ \la{218159h37}
          &=&I_n+z\E D\E Q+zp^{-1}(M_1
         + p\E \oo{D} Q).
          \eeqy   
   Now, as clearly $\|M_1\|=O(1)$ and, by Lemma \re{247153}, 
   $$\|\E \oo{D} Q\|\le \E \|\oo{D}\|\| Q\|\le \ff{|\Im z|}\E \|\oo{D}\|\le \ff{|\Im z|} \sum_{a=1}^k\sqrt{\op{Var}(p^{-1}\tr (\tQ C_a) )}=O(p^{-1})$$
       it follows that 
       $$-z\E Q = I_n+z\E D\E Q+O_{\on}(p^{-1}),$$ and
              $$-z\E Q(I_n+\E D) = I_n+O_{\on}(p^{-1}),$$
	      which implies
              $$\E Q=-z^{-1}(I_n+\E D)^{-1}+O_{\on}(p^{-1}).$$  This proves \eqre{181515h}.

        \subsubsection{Computations on $\tQ$}
 
  By the resolvent identity, we have here $\tQ=z^{-1}WW^\trans \tQ-z^{-1}I_p$ so that, by Stein Lemma again,  
 \beq \E \tQ_{ii'}&=&z^{-1}  (WW^\trans \tQ)_{ii'}-\del_{ii'}z^{-1}\\ &=&-\del_{ii'}z^{-1}+z^{-1}\E\sum_{j,l} w_{ij}w_{lj}\tQ_{li'}\\
 &=&-\del_{ii'}z^{-1}+z^{-1}\sum_{j,l,m}\E w_{ij}w_{mj} \E \pa_{w_{mj}}w_{lj}\tQ_{li'}\\
  &=&-\del_{ii'}z^{-1}+(pz)^{-1}\sum_{j,l,m}C(j)_{im} \{ \one_{(m,j)=(l,j)}\E  \tQ_{li'}-\E w_{lj}(\tQ([\del_{rm}w_{sj}+\del_{sm}w_{rj}]_{r,s=1}^p)\tQ)_{li'}\}
 \eeq
 where we used the fact that $\E w_{ij}w_{mj} =p^{-1}C(j)_{im}$ and that $$\pa_{w_{mj}}\tQ=-\tQ(\pa_{w_{mj}} WW^\trans )\tQ=-\tQ(E_{mj}W^\trans +WE_{jm})\tQ=-\tQ([\del_{rm}w_{sj}+\del_{sm}w_{rj}]_{r,s=1}^p)\tQ.$$ 
We get, somewhat similarly as previously  \beq \E \tQ_{ii'}
  &=&-\del_{ii'}z^{-1}+(pz)^{-1}\sum_j\E (C(j)\tQ)_{ii'}-(pz)^{-1}\sum_{j,l,m,s}C(j)_{im}\E w_{lj}\tQ_{lm}w_{sj}\tQ_{si'}\\&&-(pz)^{-1}\sum_{j,l,m,r}C(j)_{im}\E w_{lj}\tQ_{lr}w_{rj}\tQ_{mi'}
  \\    &=&-\del_{ii'}z^{-1}+(pz)^{-1}\sum_j\E (C(j)\tQ)_{ii'}-(pz)^{-1}\sum_{a=1}^k \E (C_a\tQ W\mc{D}(j_a)W^\trans \tQ)_{ii'}
\\
  && -(pz)^{-1}\sum_{a=1}^k\E (C_a\tQ)_{ii'}  \tr (W\mc{D}(j_a)W^\trans \tQ)  \\
   &=&-\del_{ii'}z^{-1}+z^{-1}\sum_{a=1}^k\f{c_a}{c_0}\E (C_a\tQ)_{ii'}-(pz)^{-1}\sum_{a=1}^k \E (C_a\tQ W\mc{D}(j_a)W^\trans \tQ)_{ii'}   \\&& -(pz)^{-1}\sum_{a=1}^k\E (C_a\tQ)_{ii'}  \tr (\mc{D}(j_a)W^\trans \tQ W)     
         \eeq
   Using $W^\trans \tQ W=W^\trans WQ= zQ+I_n$, we then obtain 
    \beq \E \tQ_{ii'}   &=&-\del_{ii'}z^{-1}+z^{-1}\sum_{a=1}^k\f{c_a}{c_0}\E (C_a\tQ)_{ii'}-(pz)^{-1}\sum_{a=1}^k \E(C_a\tQ W\mc{D}(j_a)W^\trans \tQ)_{ii'}   \\&& -\sum_{a=1}^k\E (C_a\tQ)_{ii'}  \ff{p}\tr (\mc{D}(j_a)Q)      -z^{-1}\sum_{a=1}^k\f{c_a}{c_0}\E (C_a\tQ)_{ii'}  \\  &=&-\del_{ii'}z^{-1} -(pz)^{-1}\sum_{a=1}^k \E(C_a\tQ W\mc{D}(j_a)W^\trans \tQ)_{ii'}     -\sum_{a=1}^k\E (C_a\tQ)_{ii'}  \ff{p}\tr (\mc{D}(j_a)Q).      
           \eeq    
           Thus for \be\la{def:M2} M_2:=\ff{z}\sum_{a=1}^k \E(C_a\tQ W\mc{D}(j_a)W^\trans \tQ),  
           \ee
           we have $$-z\E \tQ=I+z\sum_{a=1}^k\E \ff{p}\tr (\D_aQ)C_a\tQ+zp^{-1}M_2,$$\ie 
           $$\E(I_p+\sum_{a=1}^p\ff{p}\tr (Q\D_a)C_a )\tQ=-z^{-1}(I_p-p^{-1}M_2),$$ so that, for $$ {K}:=\E \ff{p}\tr (\D_aQ) \tQ-\E \ff{p}\tr (\D_aQ)\E\tQ,$$ we have 
           \beqy\nonumber \E \tQ&=&(\E(I_p+\sum_{a=1}^p\ff{p}\tr (Q\D_a)C_a ))^{-1}\E(I_p+\sum_{a=1}^p\ff{p}\tr (Q\D_a)C_a )\E \tQ\\
           \nonumber&=&(\E(I_p+\sum_{a=1}^p\ff{p}\tr (Q\D_a)C_a ))^{-1}\lf(\E(I_p+\sum_{a=1}^p\ff{p}\tr (Q\D_a)C_a )\tQ- {K}\ri)\\
           \nonumber &=& -(\E(zI_p+\sum_{a=1}^p\f{z}{p}\tr (Q\D_a)C_a ))^{-1}+p^{-1}(\E(zI_p+\sum_{a=1}^p\f{z}{p}\tr (Q\D_a)C_a ))^{-1}M_2\\
           &&-(\E(zI_p+\sum_{a=1}^p\f{z}{p}\tr (Q\D_a)C_a ))^{-1} {K}
           \la{23081515h}\eeqy 
	   Note that for any $a\in\{1,\ldots,k\}$, by Lemma \re{247153}, $$ \|{K} 
           \|=\lf\|\E\lf[  \lf(\ff{p}\tr (\D_aQ) -\E\ff{p}\tr (\D_aQ)\ri)\tQ\ri]\ri\|  \le  \ff{|\Im z|}\sqrt{\op{Var}\ff{p}\tr (\D_aQ)}=O(p^{-1}).
$$
Besides, we also have   $\|M_2\|=O(1)$ and $$\left\|\left(\E \left(zI_p+\sum_{a=1}^p zp^{-1}\tr (Q\D_a)C_a \right)\right)^{-1}\right\|\le |\Im z|^{-1},$$
          hence  
                 \beq \E \tQ  
        &=& -z^{-1}\left(I_p+\sum_{a=1}^k\E \ff{p}\tr (\mc{D}(j_a)Q)C_a\right) ^{-1}  +O_{\on}(p^{-1}).
                \eeq
This proves \eqre{181515h1}.
\subsubsection{Consequences of the loop equations: proof of Lemma \re{lem:2881516h}}\la{sec:consloopeq} We have proved  \beq \E Q&=&-z^{-1}(I_n+\E D)^{-1}+O_{\on}(p^{-1})\\
  \E \tQ  
        &=& -z^{-1}\left(I_p+\sum_{a=1}^k\E \ff{p}\tr (\mc{D}(j_a)Q)C_a\right) ^{-1}  +O_{\on}(p^{-1})
\eeq
for $D=\sum_{a=1}^k \ff{p}\tr (\tQ C_a)\D_a$. This is precisely the content of Lemma~\re{lem:2881516h}. It also implies that, for the key complex quantities  \be\la{238151}G_a=G_a(z):=\E\ff{c_0n_a}\tr (\mc{D}(j_a)Q)=\E\ff{pc_a}\tr (\mc{D}(j_a)Q),\quad \tG_a(z):=\E\ff{p}\tr(\tQ C_a),\ee we have 
		\begin{align}
			c_0G_a&=-z^{-1}(1+\tG_a)^{-1}+O(p^{-1}) \nonumber\\
			\tG_a &= -z^{-1}\ff{p}\tr C_a(I+\sum_{b=1}^kc_bG_bC_b)^{-1} +O(p^{-1}) \label{eq:tGa},
		\end{align}
     which implies that \beqy\la{208151} c_0G_a&=&\ff{-z+\ff{p}\tr C_a(I+\sum_{b=1}^kc_bG_bC_b)^{-1}  +O(p^{-1})}.
     \eeqy

 \subsection{Proof  of Proposition~\re{prop:288151}} It follows from Lemma \re{lem:keyapproxlemma} and from \eqre{208151} that, for $G_a$ as defined in \eqre{238151} and for     $z\in \C\bck\R$, \be\la{ga-Ga}|g_a(z)-G_a(z)|=O(p^{-1}).\ee  
 From Lemma \re{lem:2881516h}, we deduce that, with the notations of Proposition~\re{prop:288151}, \be\la{2881520h}\E Q_z- \bar{Q}_z =O_{\on}(p^{-1}),\qquad \E\tilde{Q}_z- \bar{\tilde{Q}}_z=O_{\on}(p^{-1}).\ee By the concentration Lemma \ref{247153}, we immediately deduce Proposition~\re{prop:288151} as long as $z$ stays at a macroscopic distance from $\R$.

 To extend the result to all $z$'s taken at a macroscopic distance from $\mc{S}\cup \{0\}$, we shall prove next that the spectrum of $W^\trans W$ remains almost surely away from $\mc{S}\cup\{0\}$ (thus proving in passing \eqref{181513h29FOC} of Theorem~\ref{Th:mu}). Let $\mathcal I$ be a closed interval of $\R$ at a macroscopic distance from $\mc{S}\cup\{0\}$. There is $\eps>0$ \st the distance from $\mathcal I$ to $\mc{S}\cup\{0\}$ is at least $2\eps$. Let $\eta>0$ \st $(2\eps^2+\eta^2)^{1/2}-\eta=\eps$. We have $$\sup_{x\in\mathcal I}\| \bar{Q}_{x+\ii\eta}\|\le (4\eps^2+\eta^2)^{-1/2}.$$ We deduce that for $p$ large enough, $$\sup_{x\in\mathcal I}\| \E Q_{x+\ii\eta}\|\le (3\eps^2+\eta^2)^{-1/2}.$$ Hence by measure concentration (using the arguments of the proof of \cite[Cor.~6]{KarginPTRF12}), with \pro tending to one, $$\sup_{x\in\mathcal I}\sup_{\lam\in \op{Spec}(W^\trans W)}|\lam-(x+\ii\eta)|^{-1}=\sup_{x\in\mathcal I}\| Q_{x+\ii\eta}\|\le (2\eps^2+\eta^2)^{-1/2},$$\ie that $$\inf_{x\in\mathcal I}\inf_{\lam\in \op{Spec}(W^\trans W)}|\lam-(x+\ii\eta)| \ge (2\eps^2+\eta^2)^{1/2},$$
which implies finally that  $$\inf_{x\in\mathcal I}\inf_{\lam\in \op{Spec}(W^\trans W)}|\lam-x| \ge (2\eps^2+\eta^2)^{1/2}-\eta=\eps.$$
This being true for any such interval $\mathcal I$, by the union bound and Lemma~\re{lem:boundWW}, we have the sought for result.

 \subsection{Proof of Proposition~\ref{prop:deteq2}}

 Our first interest is on $Q_{z_1}\mathcal D_aQ_{z_2}$. By the resolvent identity $W^\trans WQ_z-zQ_z=I_p$ applied to either of the two matrices $Q_{z_1}$ or $Q_{z_2}$, along with Stein's lemma and the results from Proposition~\ref{prop:288151}, we then get (technical details, similar to previous derivations, are omitted)
\begin{align*}
	\E [Q_{z_1}\mathcal D_a Q_{z_2}]_{ij} &= -\frac1{z_1}  \left[\mathcal D_a \bar{Q}_{z_2}\right]_{ij} - \sum_{b=1}^k \tilde{g}_b(z_1) \left[\mathcal D_b\bar{Q}_{z_1}\mathcal D_a \bar{Q}_{z_2}\right]_{ij} \\
	&-\frac1{z_1} \sum_{b=1}^k r_{ab}(z_1,z_2) \E \left[\mathcal D_b\bar{Q}_{z_2}\right]_{ij} +O(p^{-1})
\end{align*}
where we defined
\begin{align*}
	r_{ab}(z_1,z_2) &:= \E\left[ p^{-1}\tr \left( \mathcal D_a W^\trans \tilde{Q}_{z_1} C_b \tilde{Q}_{z_2}W \right) \right].
\end{align*}

Similarly, we find
\begin{align}
	\label{eq:tQWPDPWtQ}
	\E [\tilde{Q}_{z_1}W\mathcal D_a W^\trans \tilde{Q}_{z_2} ]_{ij} &= z_1 z_2 g_a(z_1) g_a(z_2) c_a c_0 \E [\tilde{Q}_{z_1}C_a\tilde{Q}_{z_2}]_{ij}+ O(p^{-1})
\end{align}
which introduces the term $\E[\tilde{Q}_{z_1}C_a\tilde{Q}_{z_2}]_{ij}$. This term is also similarly treated and gives
\begin{align*}
	\E\left[ [\tilde{Q}_{z_1}C_a\tilde{Q}_{z_2}]_{ij} \right] &= \left[ \bar{\tilde{Q}}_{z_1} C_a\bar{\tilde{Q}}_{z_2} \right]_{ij} + \sum_{b=1}^k r_{ba}(z_1,z_2) \left[ \bar{\tilde{Q}}_{z_1}C_b\bar{\tilde{Q}}_{z_2} \right]_{ij} + O(p^{-1}).
\end{align*}

To wrap up the various results, we need to identify precisely $r_{ab}(z_1,z_2)$. To this end, from \eqref{eq:tQWPDPWtQ}, we find
\begin{align*}
	r_{ab}(z_1,z_2) &= z_1 z_2 g_a(z_1)g_a(z_2) c_a c_0 \E \left[ \frac1p\tr \left(C_a \tilde{Q}_{z_1} C_b \tilde{Q}_{z_2} \right)\right] + O(p^{-1}) \\
	&= z_1 z_2 g_a(z_1)g_a(z_2) c_a c_0 \frac1p\tr \left( C_a \bar{\tilde{Q}}_{z_1}C_b \bar{\tilde{Q}}_{z_2} \right) \\
	&+ z_1 z_2 g_a(z_1)g_a(z_2) c_a c_0 \sum_{d=1}^k r_{da}(z_2,z_1) \frac1p\tr \left( C_d \bar{\tilde{Q}}_{z_1}C_b \bar{\tilde{Q}}_{z_2} \right) + O(p^{-1}).
\end{align*}
From the definition of $r_{ab}(z_1z_2)$, it is clear that $r_{da}(z_2,z_1)c_ag_a(z_1)g_a(z_2)=r_{ad}(z_1,z_2)c_dg_d(z_1)g_d(z_2)$. Thus, the above formula can be rewritten
\begin{align*}
	&r_{ab}(z_1,z_2) - z_1z_2c_0 \sum_{d=1}^k r_{ad}(z_1,z_2) c_dg_d(z_1)g_d(z_2) \frac1p\tr \left( C_d \bar{\tilde{Q}}_{z_1}C_b \bar{\tilde{Q}}_{z_2} \right) \\
	&= z_1z_2c_0 c_a g_a(z_1)g_a(z_2) \frac1p\tr \left( C_a \bar{\tilde{Q}}_{z_1}C_b \bar{\tilde{Q}}_{z_2} \right) + O(p^{-1}).
\end{align*}
This can be further rewritten under a matrix form which, after basic manipulations, leads finally to
\begin{align*}
	r_{ab}(z_1,z_2) &= [R(z_1,z_2)]_{ab} + O(p^{-1})
\end{align*}
with $R(z_1,z_2)$ defined in the statement of the lemma. The fact that $(I_k-\Omega(z_1,z_2))^{-1}$ in the expression of $R(z_1,z_2)$ is well-defined as an inverse matrix is a consequence of $\rho(\Omega(z_1,z_2))<1$ by Lemma~\ref{lem:keyapproxlemma}-(ii) with $\varepsilon=\eta=0$, for every $z_1,z_2\in\C\setminus\R$. The proof of the proposition is then completed by applying the concentration result from Lemma~\ref{247153}.

\section{Proof of Theorem \re{Th:mu}}
 
The proof of \eqref{eq:weak_conv} follows directly from Proposition~\re{prop:288151}. As for the proof of \eqref{181513h29FOC}, it was already obtained in the proof of Proposition~\ref{prop:288151}. Let us then prove the remaining second part of Theorem~\re{Th:mu}.

For $z\in \C^+$, denote $g(z)=(g_1(z), \ld, g_k(z))^\trans$ with $g_1(z), \ld, g_k(z)$ defined by Proposition~\ref{defg_as}. By Lemma \re{lem:keyapproxlemma}-(i) with $\varepsilon=0$, $\eta=0$, we have, for any $z_1,z_2\in \C\bck\R$, \begin{align}\label{eq:Gamma_g1_g2}
	\left( I_k - \Omega(z_1,z_2) \right)(g(z_1)-g(z_2)) &= (z_1-z_2)c_0 g(z_1)\odot g(z_2)
\end{align}
 where  \begin{align*}
	\Omega(z_1,z_2) &:=c_0 z_1z_2 \left\{ c_b g_a(z_1)g_a(z_2) \frac1p\tr C_a \bar{\tilde{Q}}_{z_1}C_b\bar{\tilde{Q}}_{z_2} \right\}_{a,b=1}^k.
\end{align*}

However, it is not convenient for our present investigation to work with $\Omega(z_1,z_2)$ which does not exhibit enough symmetry. We shall then proceed next by left-multiplying both sides of \eqref{eq:Gamma_g1_g2} by $\diag(c)^{\frac12}\diag(g(z_1)\odot g(z_2))^{-\frac12}$, where   the complex square root is defined thanks to the natural definition of the argument on $\C\bck\R_+$ (resp. $\C\bck\R_-$) if $\Im z_1\Im z_2>0$ (resp. if $\Im z_1\Im z_2<0$). This entails
\begin{align*}
	(I_k-\Upsilon(z_1,z_2)) \left\{ \sqrt{c_a} \frac{g_a(z_1)-g_a(z_2)}{\sqrt{g_a(z_1)g_a(z_2)}} \right\}_{a=1}^k = (z_1-z_2) c_0 \left\{\sqrt{c_ag_a(z_1)g_a(z_2)} \right\}_{a=1}^k
\end{align*}
where we defined
\begin{align*}
	\Upsilon(z_1,z_2)_{ab} &:= c_0 z_1z_2 \sqrt{c_ac_b} \sqrt{g_a(z_1)g_b(z_1) g_a(z_2) g_b(z_2)} \frac1p\tr C_a \bar{\tilde{Q}}_{z_1}C_b\bar{\tilde{Q}}_{z_2}.
\end{align*}
The matrix $\Upsilon(z_1,z_2)$ is ``more'' symmetrical than $\Omega(z_1,z_2)$ but satisfies only $\Upsilon(z_1,z_2)_{ab}=\Upsilon(z_2,z_1)_{ba}$, which shall not be good enough in what follows. To symmetrize this expression further, observe that, exchanging $z_1$ and $z_2$, we also get
\begin{align*}
	(I_k-\Upsilon(z_2,z_1)) \left\{ \sqrt{c_a} \frac{g_a(z_1)-g_a(z_2)}{\sqrt{g_a(z_1)g_a(z_2)}} \right\}_{a=1}^k = (z_1-z_2) c_0 \left\{ \sqrt{c_a g_a(z_1)g_a(z_2)} \right\}_{a=1}^k
\end{align*}
so that, summing up the two equations leads to
\begin{align*}
	(I_k-\Xi(z_1,z_2)) \left\{ \sqrt{c_a} \frac{g_a(z_1)-g_a(z_2)}{\sqrt{g_a(z_1)g_a(z_2)}} \right\}_{a=1}^k = (z_1-z_2) c_0 \left\{ \sqrt{c_a g_a(z_1)g_a(z_2)} \right\}_{a=1}^k
\end{align*}
where
\begin{align*}
	\Xi(z_1,z_2) &:= \frac12 \left( \Upsilon(z_1,z_2) +  \Upsilon(z_2,z_1) \right).
\end{align*}
In particular,
\begin{align*}
	\left( I_k - \Xi(z,z^*) \right) \left\{ \sqrt{c_a} \frac{\Im(g_a(z))}{|g_a(z)|} \right\}_{a=1}^k = \Im (z)c_0 \left\{ \sqrt{c_a} |g_a(z)| \right\}_{a=1}^k
\end{align*}
where $\Xi(z,z^*)$ is real positive and symmetric. Hence, by Lemma~\ref{lem:nonnegative_matrices}, we may take $x$ with positive entries a left eigenvector of $\Xi(z,z^*)$ with eigenvalue $\rho(\Xi(z,z^*))$. Multiplying by $x$ on the left, we get $\rho(\Xi(z,z^*))<1$. Thus, $\Xi(z,z^*)$ is invertible for every $z\in\C^+$ and we thus have
\begin{align}
	\label{eq:Xi_zz*}
	\left\{ \sqrt{c_a} \frac{\Im(g_a(z))}{|g_a(z)|} \right\}_{a=1}^k = \Im(z) c_0 \left( I_k - \Xi(z,z^*) \right)^{-1} \left\{ \sqrt{c_a} |g_a(z)| \right\}_{a=1}^k.
\end{align}

Using now the fact that
\begin{align*}
	\left| \tr (AB+CD) \right|^2 &= \left| \tr \left(\begin{bmatrix} A & C\end{bmatrix} \begin{bmatrix} B \\ D \end{bmatrix}\right)  \right|^2 \\
	&\leq \tr \left(\begin{bmatrix} A & C\end{bmatrix} \begin{bmatrix} A^* \\ C^* \end{bmatrix}\right) \tr \left(\begin{bmatrix} B^* & D^*\end{bmatrix} \begin{bmatrix} B \\ D \end{bmatrix}\right)  \\
	&= \tr (AA^*+CC^*) \tr (BB^*+DD^*)
\end{align*}
applied to $A=C_a^{\frac12}\bar{\tilde{Q}}_{z_1}C_b^{\frac12}$, $B=C_b^{\frac12}\bar{\tilde{Q}}_{z_2}C_a^{\frac12}$, $C=C_b^{\frac12}\bar{\tilde{Q}}_{z_1}C_a^{\frac12}$, and $D=C_a^{\frac12}\bar{\tilde{Q}}_{z_2}C_b^{\frac12}$, 
we find that
\begin{align*}
	|\Xi(z_1,z_2)_{ab}|^2 &\leq \Xi(z_1,z_1^*)_{ab} \Xi(z_2,z_2^*)_{ab}
\end{align*}
and thus, from Lemma~\ref{lem:Cauchy_Schwarz_spectral_radius}, we get that $\rho(\Xi(z_1,z_2))<1$ for each $z_1,z_2\in\C^+$. But since $\Xi(z,z^*)_{ab}\leq \|\Xi(z,z^*)\|=\rho(\Xi(z,z^*))<1$ for symmetric matrices, we have in addition $|\Xi(z_1,z_2)_{ab}|^2\leq 1$ for each $a,b$ so that, by e.g., $A^{-1}=\frac{ {\rm adj}(A)}{\det A}$, we finally get
\begin{align}
	\label{eq:spectralnorm_versus_radius}
	\| (I_k - \Xi(z_1,z_2))^{-1}\| &\leq \frac{K}{ |1-\rho(\Xi(z_1,z_2))|^k } 
\end{align}
for some constant $K>0$, and in particular
\begin{align}
	\label{eq:Xi_g1_g2}	
	\left\{ \sqrt{c_a} \frac{g_a(z_1)-g_a(z_2)}{\sqrt{g_a(z_1)g_a(z_2)}} \right\}_{a=1}^k = (z_1-z_2) c_0 (I_k-\Xi(z_1,z_2))^{-1} \left\{ \sqrt{c_a g_a(z_1)g_a(z_2)} \right\}_{a=1}^k.
\end{align}

\medskip

With this identity at hand, we shall show that $g(z)$ admits a limit as $z\in\C^+\to x\in\R^*$. This will be sufficient by \cite[Theorems~2.1--2.2]{CHO95} to ensure that $\mu$ admits a continuous density on $\R^*$.

Recall first the notation $$\tg_a(z):=\ff{p}\tr  C_a \bar{\tilde Q}_z,\quad \bar{\tilde{Q}}_z = -z^{-1}\left(I_p + \sum_{b=1}^k c_b g_b(z) C_b\right)^{-1}$$ (so that $c_0g_a(z)=-z^{-1}(1+\tg_a(z))^{-1}$). Then we have the following first result. 
\beg{lem}For any $\eps>0$, $g(z)$ is bounded on $\{z\in \C^+\ste |z|>\eps\}$. \en{lem}
\bpr
Note first that, by the inequality $|\tr AB^*|^2\leq \tr AA^* \tr BB^*$ with $B=I$,
\begin{align*}
	|\tilde{g}_a(z)|^2 &= \left|\frac1p\tr C_a^{\frac12} \bar{\tilde{Q}}_z C_a^{\frac12}\right|^2 \leq \frac1p\tr C_a\bar{\tilde{Q}}_zC_a\bar{\tilde{Q}}_{z^*}
\end{align*}
so that
\begin{align*}
	c_0c_a \left|z g_a(z)\tilde{g}_a(z) \right|^2 &\leq c_0 c_a |zg_a(z)|^2 \frac1p\tr C_a\bar{\tilde{Q}}_zC_a\bar{\tilde{Q}}_{z^*} = \Xi(z,z^*)_{aa}.
\end{align*}
Since $\Xi(z,z^*)_{aa}\leq \rho(\Xi(z,z^*))<1$, we thus get that for each $z\in\C^+$, $c_0 c_a |z g_a(z)\tilde{g}_a(z)|^2<1$. Hence, if $|g_a(z_n)|\to\infty$ on some sequence with $|z_n|>\eps$, this implies that $|\tilde{g}_a(z_n)|\to 0$. But by definition, $|g_a(z_n)|=|z_nc_0(1+\tilde{g}_a(z_n))|^{-1}$, which is thus bounded, contradicting the assumption. We conclude that $g_a(z)$ must remain bounded on $\{z\in \C^+\ste |z|>\eps\}$.\epr

\beg{lem} Under the additional assumptions of Theorem~\ref{Th:mu}, for any $x_0\in \R^*$, $g(z)$ admits a finite limit as $z\in \C^+$ tends to $x_0$.
 \en{lem}

\bpr  If not, by the previous lemma, one can find two sequences   $z^1_n,z^2_n \in\C^+$ tending to $ x_0\in \R^*$ \st  $g(z^1_n)\to g^1$ and $g(z_2^n)\to g^2$, $g^1\neq g^2$. From \eqref{eq:Xi_g1_g2},
\begin{align}
	\label{eq:Xi_z1n_z2n}
	\left\{ \sqrt{c_a} \frac{g_a(z^1_n) - g(z^2_n)}{\sqrt{g_a(z_1^n)g_a(z^2_n)}} \right\}_{a=1}^k &=  (z^1_n-z^2_n) c_0\left( I_k - \Xi(z^1_n,z^2_n) \right)^{-1} \left\{ \sqrt{c_a g_a(z^1_n)g_a(z^2_n)}\right\}_{a=1}^k.
\end{align}
Since $z^1_n,z^2_n \to x_0$, $z^1_n-z^2_n\to 0$. Also, since $g(z^1_n),g(z^2_n)$ are bounded by the previous lemma, we get that $(z^1_n-z^2_n) c_0\sqrt{c_a g_a(z^1_n)g_a(z^2_n)}\to 0$. It thus remains to show that $\left( I_k - \Xi(z^1_n,z^2_n) \right)^{-1}$ has uniformly bounded spectral norm, which, by \eqre{eq:spectralnorm_versus_radius},  is equivalent to showing that $$\limsup_n\rho(\Xi(z^1_n,z^2_n))<1.$$

Recall first that we obtained, from Lemma~\ref{lem:Cauchy_Schwarz_spectral_radius} and $$|\Xi(z^1_n,z^2_n)_{ab}|^2\leq \Xi(z^1_n,(z^1_n)^*)_{ab}\Xi(z^2_n,(z^2_n)^*)_{ab},$$ that $\rho(\Xi(z^1_n,z^2_n))\leq \sqrt{\rho(\Xi(z^1_n,(z^1_n)^*))\rho(\Xi(z^2_n,(z^2_n)^*))}$. Since $\rho(\Xi(z,z^*))<1$ for each $z\in\C^+$, in the limit, this only ensures that $\limsup_n\rho(\Xi(z^1_n,z^2_n))\leq 1$. We may thus show that the inequality $|\Xi(z^1_n,z^2_n)_{aa}|^2\leq \Xi(z^1_n,(z^1_n)^*)_{aa}\Xi(z^2_n,(z^2_n)^*)_{aa}$ is strict, uniformly in $n$,  for each $a$. To this end, we shall use the second part of Lemma~\ref{lem:Cauchy_Schwarz_spectral_radius}.

Letting $U^1_n:=C_a^{\frac12}\bar{\tilde{Q}}_{z_n^1}C_a^{\frac12}$ and $U^2_n:=C_a^{\frac12}\bar{\tilde{Q}}_{(z_n^2)^*}C_a^{\frac12}$, we wish to show that, uniformly on $\lambda\in\C$, $$\liminf_n \tr (U^1_n-\lambda U^2_n)(U^1_n-\lambda U^2_n)^*> 0.$$ 
For this, note that, for each $\lambda\in\C$,
\begin{align*}
	\tr ( U^1_n - \lambda U^2_n )(U^1_n - \lambda U^2_n )^* &= \tr C_a^{\frac12} \bar{\tilde{Q}}_{z_n^1} \Delta^Q_n \bar{\tilde{Q}}^*_{z_n^2} C_a \bar{\tilde{Q}}_{z_n^2} (\Delta^Q_n)^* \bar{\tilde{Q}}^*_{z_n^1} C_a^{\frac12}
\end{align*}
with
\begin{align*}
	\Delta^Q_n &:= \left( \lambda z_n^1 - (z_n^2)^* \right) I_p + \sum_{i=1}^k c_i \left( \lambda z_n^1g_i(z_n^1) - (z_n^2)^*g_i( (z_n^2)^* ) \right) C_i.
\end{align*}
From the fact that $\tr ABA^*\geq \lambda_{\rm min}(B) \tr AA^*$ when $B$ is nonnegative definite, we then get
\begin{align*}
	\tr ( U^1_n - \lambda U^2_n )(U^1_n - \lambda U^2_n )^* \geq \lambda_{\rm min} \left( \bar{\tilde{Q}}^*_{z_n^2} C_a \bar{\tilde{Q}}_{z_n^2} \right) \lambda_{\rm min} \left( \bar{\tilde{Q}}_{z_n^1}^* C_a \bar{\tilde{Q}}_{z_n^1} \right) \tr \Delta^Q_n(\Delta^Q_n)^*.
\end{align*}

Exploiting the invertibility of $C_a$    along with the fact that $\|\bar{\tilde{Q}}_{z}^{-1}\|$ is bounded uniformly on $z\in\C^+$ away from zero (by the previous lemma), we then get that
\begin{align*}
	\liminf_n \lambda_{\rm min} \left( \bar{\tilde{Q}}^*_{z_n^2} C_a \bar{\tilde{Q}}_{z_n^2} \right) \lambda_{\rm min} \left( \bar{\tilde{Q}}_{z_n^1}^* C_a \bar{\tilde{Q}}_{z_n^1} \right) > 0.
\end{align*}
By the boundedness of $g $ away from zero, we also have
\begin{align*}
	\lim_n \tr \Delta^Q_n(\Delta^Q_n)^* &= \tr \Delta^Q(\Delta^Q)^*
\end{align*}
with
\begin{align*}
	\Delta^Q &= x_0\left[ \left( \lambda - 1 \right) I_p + \sum_{i=1}^k c_i \left( \lambda g_i^1 - (g_i^2)^* \right) C_i\right].
\end{align*}

By linear independence  of the matrices $C_1,\ldots,C_k,I_p$, the quantity above cannot be zero unless $\lambda=1$ and $g_i^1=(g_i^2)^*$ for each $i$. But $\Im(g_i^1),\Im(g_i^2)\geq 0$ so that this implies $g_i^1=g_i^2\in\R$ for each $i$. But this is forbidden by assumption, and thus
\begin{align}
	\label{eq:sup_lambda}
	 \inf_{\lambda\in\C} \liminf_n \tr  ( U^1_n - \lambda U^2_n )(U^1_n - \lambda U^2_n )^* > 0.
\end{align}
This ensures (possibly over a converging subsequence, which exists for all quantities here are bounded) that
\begin{align*}
	\lim_n \tr U^1_n(U^1_n)^*\tr U^2_n(U^2_n)^* > \lim_n \left|\tr U^1_nU^2_n\right|^2.
\end{align*}
Indeed, uniformly over $x\in\R$, \eqref{eq:sup_lambda} (applied to $\lambda=x/\sqrt{2}$ and $\ii x/\sqrt{2}$) ensures that
\begin{align*}
	\lim_n \tr U^1_n (U^1_n)^* + \frac12 x^2 \lim_n \tr U^2_n (U^2_n)^* &> \sqrt{2}x \lim_n \Re \left(\tr U^1_n(U^2_n)^* \right) \\
	\lim_n \tr U^1_n (U^1_n)^* + \frac12 x^2 \lim_n \tr U^2_n (U^2_n)^* &> \sqrt{2}x \lim_n \Im \left(\tr U^1_n(U^2_n)^* \right).
\end{align*}
Taking squares left and right on both equations, summing, and taking square-roots left and right on the result, this gives, uniformly on $x$,
\begin{align*}
	\lim_n \tr U^1_n (U^1_n)^* + \frac12 x^2 \lim_n \tr U^2_n (U^2_n)^* - \sqrt{2} x \lim_n \left| \tr U^1_n(U^2_n)^* \right| > 0,
\end{align*}
the left-hand side of which is a polynomial in $x$ with discriminant $2 \lim_n \left| \tr U^1_n(U^2_n)^* \right|^2 - 2 \lim_n \tr U^1_n (U^1_n)^*\lim_n \tr U^2_n (U^2_n)^*$ which is positive.

All this finally proves, by Lemma~\ref{lem:Cauchy_Schwarz_spectral_radius}, that $\limsup_n \rho(\Xi(z^1_n,z^2_n))<1$ and therefore, recalling \eqref{eq:spectralnorm_versus_radius}, the left-hand side of \eqref{eq:Xi_z1n_z2n} converges to zero as $n\to\infty$, and so must the left-hand side. But since $g(z)$ is bounded away from zero, this implies that $g^1=g^2$, which goes against the assumption.\epr

By \cite[Theorems 2.1--2.2]{CHO95}, we then get that $\Im(g_a(z))$ is continuous on $\R^*$ and $\nu_a$ has continuous derivative $f_a(x)=\frac1\pi\Im(g_a(x))$. As $\mu=\sum_{a=1}^k c_a\nu_a$, the result follows.

\appendix
 
\section{Multidimensional Stein formula}
 \beg{lem}\la{SteinMultidim}Let $X=(X_1, \ld, X_d)$ be a centered Gaussian vector and $f : \R^d  \to \R$ be a $\mc{C}^1$  function with   derivatives having at most polynomial growth.
 Then for all $i_0=1, \ld, d$, $$\E[X_{i_0}f(X_1, \ld, X_d)]=\sum_{k=1}^d\E[X_{i_0}X_k]  \E[(\partial_kf)(X_1, \ld, X_d)].$$
 \en{lem}
 
 \bpr
 If the covariance matrix $C$ of the $X_i$'s is $I$, then the result follows from a one-dimensional integration by parts. For a more general covariance matrix $C$, introduce a standard Gaussian vector  $Y$, so that as $X\eqlaw AY$ for $A:=C^{1/2}$ and the function  $$g(y_1, \ld, y_d):=f\circ A(y_1, \ld ,y_d).$$ Then by the $C=I$ case, we have 
 \beq \E[X_{i_0}f(X_1, \ld, X_d)]&=&\sum_{j}A_{i_0j}\E[Y_{j}g(Y_1, \ld, Y_d)]\\
 &=& \sum_{j}A_{i_0j}\E[(\pa_j g)(Y_1, \ld, Y_d)]\\
  &=& \sum_{j}A_{i_0j}\sum_kA_{kj}\E[(\pa_k f)(X_1, \ld, X_d)]\\
    &=& \sum_kC_{i_0k}\E[(\pa_k f)(X_1, \ld, X_d)]
 \eeq
 \epr

\section{Concentration}
The following lemma can be found for example in \cite[Sec. 4.4.1]{agz}.
\beg{lem}\la{247153}Let $X=(X_1, \ld, X_d)$ be a standard real Gaussian vector and $f : \R^d  \to \R$ be a $\mc{C}^1$  function with  gradient $\nabla f$. Then we have \be\la{247151}\op{Var}(f(X))\;\le\; \E \|\nabla f(X)\|^2,\ee where $\|\,\cdot\,\|$ denotes the standard Euclidian norm. 

Besides, if $f$ is $k$-Lispchitz, then for any $t>0$, we have   \be\la{247152}\p (|f(X)-\E f(X)|\ge t)\; \le \; 2\me^{-\f{t^2}{2k^2}}.
\ee
\en{lem}

To apply this lemma, we shall use the following lemma. All matrix spaces, here,  are endowed with the norm $\sqrt{\Tr MM^*}$.
\beg{lem}Let $f$ be a real (resp. complex) function on $\R^+$ \st $x\mapsto f(x^2)$ is $c$-Lipschitz. Then the functions $\vfi$, $\psi$, defined on the set of $p\ti n$ complex matrices by $\vfi(X)=f(XX^*)$ and $\psi(X)=f(X^*X)$ are $c$-Lipschitz (resp. $2c$-Lipschitz).
\en{lem}

\bpr  The complex case is directly deduced from the real one by writing $f=\Re(f)+\Im(f)$. So let us suppose that $f$ is real-valued. Let $g:x\mapsto f(x^2)$ and $N:=p+n$. 
We know, by \cite[Lem. A.2]{capJTP},  that the extension of $g$ to the set of $N\ti N$ Hermitian matrices is $c$-Lipschitz.
Then, the conclusion follows from the fact that for any $p\ti n$ complex matrix $X$, $\vfi(X)$ and $\psi(X)$ are the respective $p\ti p$ upper-left corner and $n\ti n$ lower-right corner of the $N\ti N$ matrix $g(M)$, with $$M:=\bpm 0&X\\ X^*&0\epm.$$
\epr

\section{Nonnegative matrices}
The results stated here can be found in \cite{HORNJOHNSON,HORNJOHNSON2}.

\beg{lem}[Nonnegative Matrices and Dominant Eigenvectors]
	\label{lem:nonnegative_matrices}
	If $A\in\R_+^{n\times n}$ is nonnegative (resp., positive), then $\rho(A)$ is an eigenvalue of $A$ having an eigenvector with nonnegative (resp., positive) entries.
\end{lem}

\begin{lem}[Spectral Radii]
	\label{lem:spectral_radius} 
	Let $A,B\in\R^{n\times n}$ be such that $|A_{ij}|\leq B_{ij}$ for all $1\leq i,j\leq n$. Then, with $\rho$ the spectral radius,
	\begin{align*}
		\rho(A)\leq \rho(B).
	\end{align*}
\end{lem}

\begin{lem}[Cauchy-Schwarz for Spectral Radii (adapted from {\cite[Lemma~5.7.9]{HORNJOHNSON2}})]
	\label{lem:Cauchy_Schwarz_spectral_radius}
	Let $A,B\in\R_+^{n\times n}$ be non negative matrices and $C\in\R^{n\times n}$ be such that $C_{ij}\leq \sqrt{A_{ij}B_{ij}}$. Then,
	\begin{align*}
		\rho(C)\leq \sqrt{\rho(A)\rho(B)}.
	\end{align*}
	Besides, if, for each $i$, either both the $i$-th row and the $i$-th column of $C$ are null or there exists $j$ such that $C_{ij}<\sqrt{A_{ij}B_{ij}}$, then the inequality is strict.
\end{lem}

\en{document}

%% file: Figure_Histogram.tex
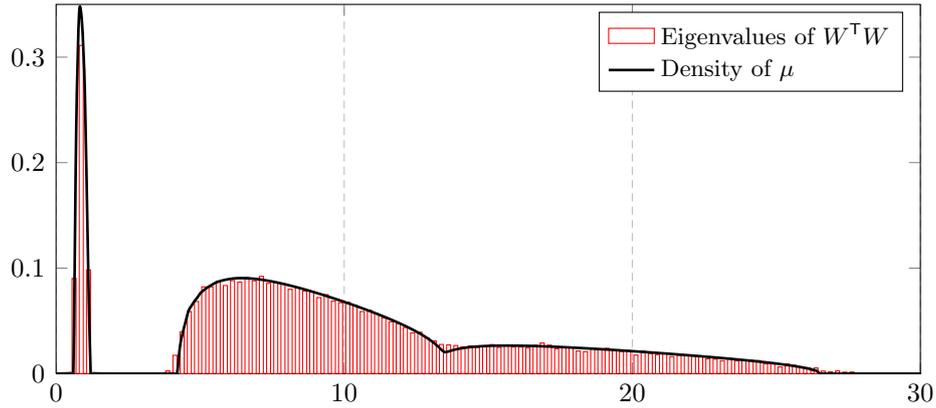
\begin{figure}[t]
  \centering
  \begin{tikzpicture}[font=\footnotesize]
    \renewcommand{\axisdefaulttryminticks}{4} 
    \tikzstyle{every major grid}+=[style=densely dashed]       
    \tikzstyle{every axis y label}+=[yshift=-10pt] 
    \tikzstyle{every axis x label}+=[yshift=5pt]
    \tikzstyle{every axis legend}+=[cells={anchor=west},fill=white,
        at={(0.98,0.98)}, anchor=north east, font=\scriptsize ]
    \begin{axis}[
      width=.8\linewidth,
      height=.4\linewidth,
      xmin=0,
      ymin=0,
      xmax=30,
      ymax=.35,
      xtick={0,10,20,30,40},
      bar width=1.5pt,
      grid=major,
      ymajorgrids=false,
      ]
      \addplot[area legend,ybar,mark=none,color=red,fill=white] coordinates{
	      (0.125000,0.000000)(0.375000,0.000000)(0.625000,0.090250)(0.875000,0.311125)(1.125000,0.098250)(1.375000,0.000375)(1.625000,0.000000)(1.875000,0.000000)(2.125000,0.000000)(2.375000,0.000000)(2.625000,0.000000)(2.875000,0.000000)(3.125000,0.000000)(3.375000,0.000000)(3.625000,0.000500)(3.875000,0.002625)(4.125000,0.017375)(4.375000,0.039625)(4.625000,0.058875)(4.875000,0.068375)(5.125000,0.082250)(5.375000,0.082375)(5.625000,0.087750)(5.875000,0.083500)(6.125000,0.088250)(6.375000,0.086750)(6.625000,0.091125)(6.875000,0.088000)(7.125000,0.092250)(7.375000,0.085500)(7.625000,0.086375)(7.875000,0.085125)(8.125000,0.080000)(8.375000,0.082000)(8.625000,0.078250)(8.875000,0.078125)(9.125000,0.072125)(9.375000,0.074875)(9.625000,0.068750)(9.875000,0.067250)(10.125000,0.067750)(10.375000,0.065000)(10.625000,0.058875)(10.875000,0.060250)(11.125000,0.056750)(11.375000,0.053000)(11.625000,0.049125)(11.875000,0.047500)(12.125000,0.043750)(12.375000,0.038625)(12.625000,0.039250)(12.875000,0.033875)(13.125000,0.031000)(13.375000,0.027625)(13.625000,0.027250)(13.875000,0.026625)(14.125000,0.025375)(14.375000,0.026250)(14.625000,0.024750)(14.875000,0.025875)(15.125000,0.027500)(15.375000,0.024875)(15.625000,0.026750)(15.875000,0.025250)(16.125000,0.025250)(16.375000,0.024500)(16.625000,0.027000)(16.875000,0.029000)(17.125000,0.027000)(17.375000,0.023625)(17.625000,0.024875)(17.875000,0.025000)(18.125000,0.021500)(18.375000,0.020875)(18.625000,0.024375)(18.875000,0.023250)(19.125000,0.024000)(19.375000,0.022250)(19.625000,0.021000)(19.875000,0.021625)(20.125000,0.017625)(20.375000,0.021500)(20.625000,0.019000)(20.875000,0.018000)(21.125000,0.019375)(21.375000,0.016125)(21.625000,0.017375)(21.875000,0.017625)(22.125000,0.016625)(22.375000,0.015875)(22.625000,0.014000)(22.875000,0.013750)(23.125000,0.012250)(23.375000,0.014375)(23.625000,0.012125)(23.875000,0.012750)(24.125000,0.012625)(24.375000,0.012375)(24.625000,0.009250)(24.875000,0.009250)(25.125000,0.006375)(25.375000,0.009125)(25.625000,0.007750)(25.875000,0.006125)(26.125000,0.005375)(26.375000,0.005375)(26.625000,0.002375)(26.875000,0.001750)(27.125000,0.002625)(27.375000,0.001375)(27.625000,0.001500)(27.875000,0.000750)(28.125000,0.000375)(28.375000,0.000375)(28.625000,0.000500)(28.875000,0.000250)(29.125000,0.000000)(29.375000,0.000000)(29.625000,0.000000)(29.875000,0.000125)(30.125000,0.000000)(30.375000,0.000000)(30.625000,0.000000)(30.875000,0.000000)(31.125000,0.000000)(31.375000,0.000000)(31.625000,0.000000)(31.875000,0.000000)(32.125000,0.000000)(32.375000,0.000000)(32.625000,0.000000)(32.875000,0.000000)(33.125000,0.000000)(33.375000,0.000000)(33.625000,0.000000)(33.875000,0.000000)(34.125000,0.000000)(34.375000,0.000000)(34.625000,0.000000)(34.875000,0.000000)(35.125000,0.000000)
      };
      \addplot[black,smooth,line width=1pt] coordinates{
	      (0.,0.000058)(0.1,0.000074)(0.2,0.000099)(0.3,0.000139)(0.4,0.000214)(0.5,0.000388)(0.6,0.001071)(0.7,0.232802)(0.8,0.341812)(0.9,0.332563)(1.,0.266492)(1.1,0.122009)(1.2,0.000711)(1.3,0.000305)(1.4,0.00018)(1.5,0.000123)(1.6,0.00009)(1.7,0.00007)(1.8,0.000057)(1.9,0.000048)(2.,0.000041)(2.1,0.000036)(2.2,0.000033)(2.3,0.00003)(2.4,0.000028)(2.5,0.000027)(2.6,0.000025)(2.7,0.000025)(2.8,0.000024)(2.9,0.000024)(3.,0.000025)(3.1,0.000025)(3.2,0.000026)(3.3,0.000027)(3.4,0.000029)(3.5,0.000032)(3.6,0.000035)(3.7,0.000039)(3.8,0.000045)(3.9,0.000055)(4.,0.00007)(4.1,0.000103)(4.2,0.000254)(4.3,0.025964)(4.4,0.039037)(4.5,0.049816)(4.6,0.059829)(4.7,0.064288)(4.8,0.068613)(4.9,0.071997)(5.,0.075977)(5.1,0.078521)(5.2,0.080684)(5.3,0.082496)(5.4,0.084063)(5.5,0.08566)(5.6,0.08686)(5.7,0.087764)(5.8,0.088496)(5.9,0.089087)(6.,0.089581)(6.1,0.090011)(6.2,0.090313)(6.3,0.090469)(6.4,0.090522)(6.5,0.090497)(6.6,0.090409)(6.7,0.09027)(6.8,0.09008)(6.9,0.08983)(7.,0.089517)(7.1,0.08915)(7.2,0.088739)(7.3,0.088292)(7.4,0.087812)(7.5,0.087303)(7.6,0.086763)(7.7,0.086192)(7.8,0.085591)(7.9,0.084963)(8.,0.084312)(8.1,0.083641)(8.2,0.08295)(8.3,0.082241)(8.4,0.081514)(8.5,0.08077)(8.6,0.080008)(8.7,0.079232)(8.8,0.078441)(8.9,0.077636)(9.,0.076819)(9.1,0.075989)(9.2,0.075148)(9.3,0.074294)(9.4,0.073429)(9.5,0.072553)(9.6,0.071665)(9.7,0.070768)(9.8,0.06986)(9.9,0.068941)(10.,0.068013)(10.1,0.067073)(10.2,0.066123)(10.3,0.065162)(10.4,0.064191)(10.5,0.063208)(10.6,0.062214)(10.7,0.061207)(10.8,0.060189)(10.9,0.059157)(11.,0.058112)(11.1,0.057053)(11.2,0.055978)(11.3,0.054888)(11.4,0.05378)(11.5,0.052654)(11.6,0.051508)(11.7,0.05034)(11.8,0.049148)(11.9,0.047931)(12.,0.046685)(12.1,0.045408)(12.2,0.044095)(12.3,0.042743)(12.4,0.041345)(12.5,0.039896)(12.6,0.038386)(12.7,0.036806)(12.8,0.035141)(12.9,0.033374)(13.,0.031479)(13.1,0.029421)(13.2,0.027151)(13.3,0.024612)(13.4,0.021897)(13.5,0.020354)(13.6,0.020843)(13.7,0.02169)(13.8,0.022459)(13.9,0.023113)(14.,0.023669)(14.1,0.024141)(14.2,0.024544)(14.3,0.02489)(14.4,0.025187)(14.5,0.025442)(14.6,0.025661)(14.7,0.025848)(14.8,0.026006)(14.9,0.026139)(15.,0.026249)(15.1,0.02634)(15.2,0.026411)(15.3,0.026466)(15.4,0.026506)(15.5,0.026531)(15.6,0.026543)(15.7,0.026544)(15.8,0.026533)(15.9,0.026512)(16.,0.026482)(16.1,0.026442)(16.2,0.026394)(16.3,0.026338)(16.4,0.026275)(16.5,0.026205)(16.6,0.026129)(16.7,0.026047)(16.8,0.025958)(16.9,0.025865)(17.,0.025766)(17.1,0.025662)(17.2,0.025554)(17.3,0.025441)(17.4,0.025324)(17.5,0.025203)(17.6,0.025078)(17.7,0.02495)(17.8,0.024818)(17.9,0.024683)(18.,0.024545)(18.1,0.024404)(18.2,0.024259)(18.3,0.024112)(18.4,0.023963)(18.5,0.02381)(18.6,0.023656)(18.7,0.023499)(18.8,0.02334)(18.9,0.023179)(19.,0.023015)(19.1,0.022849)(19.2,0.02268)(19.3,0.02251)(19.4,0.022338)(19.5,0.022164)(19.6,0.021988)(19.7,0.021811)(19.8,0.021633)(19.9,0.021453)(20.,0.02127)(20.1,0.021086)(20.2,0.0209)(20.3,0.020711)(20.4,0.020521)(20.5,0.02033)(20.6,0.020137)(20.7,0.019943)(20.8,0.019748)(20.9,0.019552)(21.,0.019355)(21.1,0.019156)(21.2,0.018954)(21.3,0.018749)(21.4,0.018542)(21.5,0.018334)(21.6,0.018124)(21.7,0.017913)(21.8,0.017701)(21.9,0.017489)(22.,0.017276)(22.1,0.017062)(22.2,0.016845)(22.3,0.016622)(22.4,0.016396)(22.5,0.016168)(22.6,0.015939)(22.7,0.015708)(22.8,0.015476)(22.9,0.015242)(23.,0.015009)(23.1,0.014776)(23.2,0.014541)(23.3,0.014295)(23.4,0.01404)(23.5,0.013785)(23.6,0.013529)(23.7,0.01327)(23.8,0.013009)(23.9,0.012746)(24.,0.012483)(24.1,0.012227)(24.2,0.011955)(24.3,0.011655)(24.4,0.011361)(24.5,0.011064)(24.6,0.010763)(24.7,0.010457)(24.8,0.010148)(24.9,0.009849)(25.,0.009548)(25.1,0.009163)(25.2,0.008807)(25.3,0.008442)(25.4,0.008065)(25.5,0.007677)(25.6,0.007336)(25.7,0.006823)(25.8,0.006375)(25.9,0.005883)(26.,0.005353)(26.1,0.004925)(26.2,0.004028)(26.3,0.003341)(26.4,0.002595)(26.5,0.000205)(26.6,0.000011)(26.7,0.000008)(26.8,0.000006)(26.9,0.000005)(27.,0.000005)(27.1,0.000004)(27.2,0.000004)(27.3,0.000004)(27.4,0.000003)(27.5,0.000003)(27.6,0.000003)(27.7,0.000003)(27.8,0.000003)(27.9,0.000003)(28.,0.000003)(28.1,0.000002)(28.2,0.000002)(28.3,0.000002)(28.4,0.000002)(28.5,0.000002)(28.6,0.000002)(28.7,0.000002)(28.8,0.000002)(28.9,0.000002)(29.,0.000002)(29.1,0.000002)(29.2,0.000002)(29.3,0.000002)(29.4,0.000002)(29.5,0.000002)(29.6,0.000002)(29.7,0.000002)(29.8,0.000001)(29.9,0.000001)(30.,0.000001)(30.1,0.000001)(30.2,0.000001)(30.3,0.000001)(30.4,0.000001)(30.5,0.000001)(30.6,0.000001)(30.7,0.000001)(30.8,0.000001)(30.9,0.000001)(31.,0.000001)(31.1,0.000001)(31.2,0.000001)(31.3,0.000001)(31.4,0.000001)(31.5,0.000001)(31.6,0.000001)(31.7,0.000001)(31.8,0.000001)(31.9,0.000001)(32.,0.000001)(32.1,0.000001)(32.2,0.000001)(32.3,0.000001)(32.4,0.000001)(32.5,0.000001)(32.6,0.000001)(32.7,0.000001)(32.8,0.000001)(32.9,0.000001)(33.,0.000001)(33.1,0.000001)(33.2,0.000001)(33.3,0.000001)(33.4,0.000001)(33.5,0.000001)(33.6,0.000001)(33.7,0.000001)(33.8,0.000001)(33.9,0.000001)(34.,0.000001)(34.1,0.000001)(34.2,0.000001)(34.3,0.000001)(34.4,0.000001)(34.5,0.000001)(34.6,0.000001)(34.7,0.000001)(34.8,0.000001)(34.9,0.000001)(35.,0.000001)
      };
      \legend{ {Eigenvalues of $W^\trans W$},{Density of $\mu$} }
    \end{axis}
  \end{tikzpicture}
  \caption{Eigenvalues of $W^\trans W$ (across $1\,000$ realizations) versus $\mu$, $n=32$, $p=256$, $k=3$, $c_1=1/8$, $c_2=5/8$, $c_3=1/4$, $[C_a]_{ij}=(8(a-1)+1) [(a-1)/5]^{|i-j|}$. The density of $\mu$ is computed as detailed  in Remark \re{numerics_for_mu}.}
  \label{fig:histWW}
\end{figure}